\documentclass[%
 aip,
 amsmath,amssymb,
reprint, onecolumn 
]{revtex4-1}
\usepackage[margin=1in]{geometry}
\usepackage{graphicx} 
\usepackage{amsmath}
\usepackage{booktabs}
\usepackage{array}
\usepackage{hyperref} 
\usepackage{amssymb}
\usepackage{algpseudocode}
\usepackage{algpseudocode}
\makeatletter
\def\@email#1#2{%
 \endgroup
 \patchcmd{\titleblock@produce}
  {\frontmatter@RRAPformat}
  {\frontmatter@RRAPformat{\produce@RRAP{*#1\href{mailto:#2}{#2}}}\frontmatter@RRAPformat}
  {}{}
}%
\makeatother
\begin{document}

\preprint{AIP/123-QED}

\title[Full Vectorial Maxwell Equations with Continuous Angular Indices]{Full Vectorial Maxwell Equations with Continuous Angular Indices}
\author{Mustafa Bakr}
 \email{mustafa.bakr@physics.ox.ac.uk}
\affiliation{ 
Clarendon Laboratory, Department of Physics, University of Oxford}%

\date{\today}

\begin{abstract}
This article presents a mathematical framework for solving Maxwell's equations in spherical geometries with continuous angular indices. We extend beyond standard discrete harmonic decomposition to a continuous spectral representation using generalized spectral integrals, capturing electromagnetic solutions that exhibit singular behavior yet yield finite-energy fields at the geometric center. For continuous angular indices $\ell, m \in \mathbb{R}$, we study the existence and uniqueness of solutions in weighted Sobolev spaces $H^s_{\alpha(\ell,m)}(\Omega)$ following the framework established in ~\cite{adams2003, reed1975}, prove finite energy for $\ell > -\frac{1}{2}$, and construct explicit spectral kernels via biorthogonal function systems. The framework encompasses non-separable spherical modes where field components couple through vectorial curl operations. We present asymptotic analysis of singular field behavior, investigate convergence rates for spectral approximations, and validate the theoretical framework through Galerkin projection methods and numerical spectral integration.

\end{abstract}

\maketitle

\section{Problem Statement}
The classical analysis of electromagnetic fields in spherical coordinates relies on separable solutions using integer-order spherical harmonics \( Y_{\ell}^{m}(\theta, \phi) \)~\cite{jackson1999, stratton1941, edmonds1996, varshalovich1988, barrera1985}, assuming \( \ell \in \mathbb{Z}_+ \), \( m \in \mathbb{Z} \), and periodic azimuthal boundary conditions. While this framework has been successful for spherically symmetric problems, it fails in systems where full azimuthal periodicity is broken or where continuous symmetry breaking plays a central role. We pose a fundamental question:  What is the structure of electromagnetic fields governed by the full vectorial Maxwell equations in spherical geometry when the angular indices \( \ell \) and \( m \) are allowed to vary continuously? This leads us to a regime where the field equations become non-separable, \( \ell \) and \( m \) become coupled through curl operations, and the radial behaviour exhibits singular scaling \( E_r \sim r^{\alpha(\ell, m)} \) near the origin. To address this, We develop a generalized continuous spectral integral formalism that spans non-integer \( (\ell, m) \), enabling continuous expansions over angular momentum space. This framework: 
\begin{enumerate}
    \item Characterises the singularity structure analytically.
    \item Establishes convergence and regularity in weighted Sobolev spaces.
    \item Avoids separable \textit{ansatz} by solving the full coupled angular PDE system.
\end{enumerate}

The analysis of these singular structures in cylindrical geometry, where the variables remain separable and continuous azimuthal modes naturally arise from broken \( 2\pi \) periodicity, is presented in~\cite{bakr1, bakr2}. In this simpler case, we showed that for continuous azimuthal index \( \nu \in (0,1) \), the transverse electromagnetic field components exhibit a power-law singularity \( E_\rho, E_\phi \sim \rho^{\nu-1} \) near the axis, which remains integrable in the electromagnetic energy norm. This enables explicit harmonic decompositions of individual field components, governed by scalar Helmholtz equations with continuous-order Bessel functions. We then extend the analysis to spherical geometry, where continuous angular indices \( (\ell, m) \in \mathbb{R}^2 \) lead to coupling between components through the curl operations in Maxwell's equations. In this setting, variables are no longer separable, and a full vectorial treatment via spectral integrals over non-integer modes becomes essential.

\section{Main Contributions of This Work}
In this article, we propose a vectorial framework for Maxwell’s equations with continuous angular indices. In doing so, we move beyond classical spherical harmonics (which require integer numbers) to a continuous spectrum of modes labeled by non-integer and continuous “angular momentum” parameters. At the heart of our article is a spectral decomposition of electromagnetic fields on a spherical domain using non-integer angular indices. Mathematically, this means that instead of expanding fields in the usual discrete basis of spherical harmonics $Y_{\ell m}(\theta,\phi)$ (with $\ell\in\mathbb{N}_0$ and $-\ell\le m\le \ell$), we consider a continuous spectrum of solutions labeled by real (or complex) values of $\ell$ (and likewise treat $m$ potentially as a continuous parameter). In the context of spherical harmonics, it is known that the requirement $0\le\phi<2\pi$ (full $2\pi$ azimuthal periodicity) enforces $m$ to be an integer, and similarly regularity on $0\le\theta\le\pi$ quantises $\ell$. If one relaxes these global conditions – for example, through analytical continuation that admits multi-valued functions on appropriate Riemann surfaces – then associated Legendre functions of continuous degree can appear as legitimate solutions. 
We develop what we term a continuous integral representation of the electromagnetic field, essentially an integral over a continuum of $\ell$ and $m$ values rather than a sum. Here we show that analytical continuation of the discrete harmonic basis to continuous indices naturally yields solutions that exhibit singular behavior while maintaining finite energy under the condition $\ell > -\frac{1}{2}$.

\subsubsection{Derivations and Projection Operators}
We introduce projection operators that project the fields onto these generalized spherical harmonics with continuous indices. In classical spherical harmonic analysis, projection onto the $(\ell,m)$ mode is achieved by an integral of the field against $Y_{\ell m}^*(\theta,\phi)$ over the solid angle. For non-integer and continuous $\ell$, we must define analogous basis functions $\Phi_{\ell m}(\theta,\phi)$ and a resolution of identity. Indeed, in the Appendix D we define vector spherical functions ${\Phi}_{\ell m}$ built from associated Legendre functions $P_\ell^{|m|}(\cos\theta)$ (analytically continued to non-integer $\ell$). Orthogonality in the continuous sense is addressed by introducing a spectral weight $w(\ell,m)$ such that a generalized Parseval/Plancherel formula holds. We outline how an $L^2$ function can be expanded as $f(\theta,\phi)=\iint a(\ell,m)\Phi_{\ell m}(\theta,\phi)d\ell dm$ with $\iint |a(\ell,m)|^2 w(\ell,m)d\ell dm = |f|^2_{L^2}$. We think this is a logical extension of the familiar spherical harmonic expansion – essentially a continuous analog thereof. We carefully derive properties of these projections, showing that as one approaches the integer limits, one recovers the usual orthonormality, whereas in between one must use delta-function normalisation. While our exposition of this is quite dense, it is mathematically consistent with the theory of generalized eigenfunctions. For instance, the relationship $\ell(\ell+1)=s(1-s)$ is invoked, connecting the eigenvalue $\ell(\ell+1)$ of the spherical Laplacian to a continuous spectral parameter $s$ (this same relation arises in the theory of the non-compact hyperbolic Laplacian, where $s$ plays the role of $1/2 + i\nu$ for $\nu\in\mathbb{R}$). This correspondence allows us to import techniques from the generalized integral (which are the canonical continuous-spectrum eigenfunctions on hyperbolic manifolds) and apply them to this electromagnetic problem. We show that continuous $\ell$ solutions can be excited in geometries with modified boundary conditions, such as wedge geometries. Based on analogy to standard results (e.g. the spectral theorem for self-adjoint operators in Hilbert space), one expects that the vector Laplacian or curl-curl operator on an appropriate function space will yield a spectrum consisting of the usual discrete modes plus an integral over a continuum of improper modes. In summary, these solutions with continuous indices are not normalizable in the usual $L^2$ sense, similar to plane waves or associated Legendre functions of continuous degree. In Appendix A we address convergence using Sobolev norms to show that truncating the spectral integral yields small errors. We provide the convergence estimate $\|f - f_N\|_{L^2(S^2)} \le C \|f\|_{H^s} N^{-s+1/2+\epsilon}$ for $s>1/2$, which represents a spectral convergence rate for our approximation method. We show that the continuous-harmonic expansion converges for fields with sufficient regularity through the use of weighted Sobolev spaces and careful accounting of the singular behavior.

\subsubsection{Coupling of Vector Components}
A significant mathematical challenge in full Maxwell equations (in contrast to scalar wave equations) is the coupling introduced by the vector nature of the fields. In spherical coordinates, even in standard integer-$\ell$ analysis, one often expands the vector field in vector spherical harmonic basis functions (which mix components $E_r$, $E_\theta$, $E_\phi$). We acknowledge this complexity and derive the coupled equations for the field components by substituting the spectral ansatz into Maxwell's curl equations. We present the explicit form of the curl operator in spherical coordinates acting on these continuous-index modes. For example, one relation shown is $i\omega\mu_0 H_{\ell,m}^r = \frac{1}{r\sin\theta}\left[\frac{\partial}{\partial\theta}(E_{\ell,m}^\phi \sin\theta) - \frac{\partial}{\partial\phi}(E_{\ell,m}^\theta)\right]$ (this is analogous to the standard expression, now assumed to hold for non-integer $\ell$ as well). Such equations demonstrate how $E_r$ acts as a source for $E_\theta$ and $E_\phi$, and vice versa through angular derivatives. Importantly, we do not assume a priori that the fields factorize as $R(r)Y(\theta,\phi)$ and allow non-separable solutions. This means the radial dependence of each component can carry an implicit dependence on $\ell, m$ that is not factorizable. We introduce a generalized angular operator $L_{\text{ang}}$ which includes a potential term $V(\theta,\phi)$ arising from the coupling between different continuous-index modes. In essence, we have an eigenvalue problem $L_{\text{ang}}[{\Phi}_{\ell m}(\theta,\phi)] = \lambda(\ell,m) {\Phi}_{\ell m}$, where $\lambda$ corresponds to the separation constant (related to $\ell(\ell+1)$ in the symmetric case, now generalized). The presence of $V(\theta,\phi)$ means that $L_{\text{ang}}$'s eigenfunctions for different $(\ell,m)$ are not orthogonal in the usual way, reinforcing the need for the continuous spectrum approach.

\subsubsection{Asymptotic Behavior}
One noteworthy result is our analysis of the asymptotic behavior of these continuous modes. Near $r \to 0$, the solutions are singular. We find solutions behaving like $E_r \sim r^{\alpha(\ell,m)}$ with some exponent $\alpha$ that may be negative (indicating a singularity at the origin). Crucially, we show these singular solutions are integrable in the energy sense — that is, the electromagnetic energy density remains finite. For example, in cylindrical coordinates, we found an electric field scaling $| {E}| \sim \rho^{\nu - 1}$ near the central axis, where $\nu$ is a continuous angular parameter. In a test case we take $\nu = 0.1$, yielding $| {E}| \sim \rho^{-0.9}$. This weak singularity still leads to a finite energy (since $| {E}|^2 \sim \rho^{-1.8}$ and the area element $\rho\, d\rho$ gives an integrand $\sim \rho^{-0.8}$ near $0$, which is integrable)~\cite{bakr1}. Indeed, we emphasize that these singular electromagnetic modes are mathematically admissible — they are not divergences that violate energy conservation or blow up the field norm. The convergence criteria in weighted Sobolev spaces ensure that even though the field is singular, it lies in $H^s$ for appropriate $s < 1$ (likely $s < 1/2$ for the most singular cases), so that the energy (which corresponds to the $H^1$ norm roughly) is finite.

\section{Vectorial Maxwell Equations in Spherical Coordinates  Non-Separable Continuous Spectrum}
Having established the mathematical framework for continuous angular indices in cylindrical geometry in~\cite{bakr1, bakr2}, we now extend our analysis to spherical coordinates, where the vectorial nature of Maxwell's equations creates essential coupling between field components. This coupling destroys the separability that characterizes integer-index solutions, necessitating the full continuous spectral framework developed in this section.
\subsection{Scalar Case}
To establish the foundation for our vectorial analysis, we begin with the scalar Laplacian in spherical coordinates. For a scalar function $\Psi(r, \theta, \varphi)$ 
\begin{equation}
\nabla^2 \Psi = \frac{1}{r^2} \frac{\partial}{\partial r} \left( r^2 \frac{\partial \Psi}{\partial r} \right) + \frac{1}{r^2 \sin \theta} \frac{\partial}{\partial \theta} \left( \sin \theta \frac{\partial \Psi}{\partial \theta} \right) + \frac{1}{r^2 \sin^2 \theta} \frac{\partial^2 \Psi}{\partial \phi^2}.
\label{eq laplacian}
\end{equation}
Now, introduce the separation of variables 
\begin{equation}
\Psi(r, \theta, \varphi) = R(r) Y_{\ell}^{m}(\theta, \varphi).
\label{eq scalar_separation}
\end{equation}
The spherical harmonics are eigenfunctions of the angular momentum operators 
\begin{equation}
L^2 Y_{\ell}^{m} = \ell(\ell + 1) Y_{\ell}^{m}.
\label{eq L2_eigenvalue}
\end{equation}
\begin{equation}
L_z Y_{\ell}^{m} = m Y_{\ell}^{m}.
\label{eq Lz_eigenvalue}
\end{equation}
The Laplacian acting on \( \Psi \) becomes 
\begin{equation}
\nabla^2 \Psi = \left[ \frac{1}{r^2} \frac{d}{d r} \left( r^2 \frac{d R}{d r} \right) - \frac{\ell(\ell + 1)}{r^2} R(r) \right] Y_{\ell}^{m}(\theta, \varphi),
\label{eq laplacian_separated}
\end{equation}
where $L^2$ and $L_z$ commute with $\nabla^2$. This means the scalar Laplacian is diagonal in the basis of spherical harmonics $Y_{\ell}^{m}$.

\subsection{Vector Case}
We now consider the structure of the Laplacian acting on a vector field $ {A}(r, \theta, \varphi)$. In contrast to the scalar case, the vector Laplacian introduces coupling between components due to the curvature of the spherical coordinate basis. This coupling becomes essential when extending to continuous angular indices, as it prevents the simple factorization that characterizes integer-index solutions.

\subsubsection{Vector Laplacian in Curvilinear Coordinates}

In spherical coordinates, the Laplacian of a vector field is not the component-wise application of the scalar Laplacian. Instead, it is given by 
\begin{equation}
\nabla^2  {A} = \nabla (\nabla \cdot  {A}) - \nabla \times (\nabla \times  {A}).
\label{eq vector_laplacian_identity}
\end{equation}
This identity ensures compatibility with Maxwell's equations in vacuum. Importantly, in spherical coordinates $(r, \theta, \varphi)$, the coordinate basis vectors $\hat{ {r}}, \hat{\boldsymbol{\theta}}, \hat{\boldsymbol{\varphi}}$ vary with position, and additional connection terms emerge when applying differential operators. We write the vector field in spherical components 
\begin{equation}
 {A}(r, \theta, \varphi) = A_r(r,\theta,\varphi) \hat{ {r}} + A_\theta(r,\theta,\varphi) \hat{\boldsymbol{\theta}} + A_\varphi(r,\theta,\varphi) \hat{\boldsymbol{\varphi}}.
\label{eq vector_field_components}
\end{equation}
The divergence and curl in spherical coordinates mix these components nontrivially, even for a field initially aligned in one direction. Consequently, the Laplacian couples radial and angular components, and no single component evolves independently.

\subsubsection{Expansion in Vector Spherical Harmonics}

To analyse angular structure systematically, we expand \(  {A} \) in the vector spherical harmonics (VSH), which form a complete orthogonal basis for square-integrable vector fields on the sphere 
\begin{equation}
 {A}(r, \theta, \phi) = \sum_{\ell=0}^\infty \sum_{m=-\ell}^\ell \left[
a^{(r)}_{\ell m}(r) \,  {Y}_{\ell m}^{(r)}(\theta,\phi) +
a^{(1)}_{\ell m}(r) \,  {Y}_{\ell m}^{(1)}(\theta,\phi) +
a^{(2)}_{\ell m}(r) \,  {Y}_{\ell m}^{(2)}(\theta,\phi)
\right],
\label{eq VSH_expansion}
\end{equation}
where the three orthogonal vector harmonics are defined as 
\begin{align}
 {Y}_{\ell m}^{(r)}(\theta,\phi) &= Y_\ell^m(\theta,\phi) \, \hat{ {r}} \label{eq vsh_radial}, \\
 {Y}_{\ell m}^{(1)}(\theta,\phi) &= r \nabla Y_\ell^m(\theta,\phi) = \frac{\partial Y_\ell^m}{\partial \theta} \, \hat{\boldsymbol{\theta}} + \frac{1}{\sin \theta} \frac{\partial Y_\ell^m}{\partial \phi} \, \hat{\boldsymbol{\phi}}, \label{eq vsh_grad} \\
 {Y}_{\ell m}^{(2)}(\theta,\phi) &= \hat{ {r}} \times \nabla Y_\ell^m = \frac{1}{\sin \theta} \frac{\partial Y_\ell^m}{\partial \phi} \, \hat{\boldsymbol{\theta}} - \frac{\partial Y_\ell^m}{\partial \theta} \, \hat{\boldsymbol{\phi}}, \label{eq vsh_curl}
\end{align}
where \(  {Y}^{(r)} \) captures the longitudinal (radial) part, \(  {Y}^{(1)} \) is the polar (even-parity) component, and \(  {Y}^{(2)} \) is the axial (odd-parity) component. These satisfy orthogonality relations over the sphere and serve as the basis for studying wave equations and gauge theories on curved manifolds.

\subsubsection{Coupling under the Laplacian}

When the vector Laplacian acts on \(  {A} \), the radial and angular parts are no longer separable in general. In particular the term \( \nabla(\nabla \cdot  {A}) \) produces both radial and angular derivatives acting on \( A_\theta, A_\phi \). Similarly, the term \( \nabla \times (\nabla \times  {A}) \) introduces mixing between \( A_r \) and the transverse components through second-order angular derivatives. Thus, the Laplacian of a vector field expressed in VSH generally takes the form 
\begin{equation}
\nabla^2  {A} = \sum_{\ell,m} \left[
\mathcal{D}^{(r)}_{\ell m}(r) \,  {Y}_{\ell m}^{(r)} +
\mathcal{D}^{(1)}_{\ell m}(r) \,  {Y}_{\ell m}^{(1)} +
\mathcal{D}^{(2)}_{\ell m}(r) \,  {Y}_{\ell m}^{(2)}
\right],
\label{eq vector_laplacian_expansion}
\end{equation}
where the differential operators \( \mathcal{D}^{(i)}_{\ell m}(r) \) involve second-order radial and angular derivatives and cannot be decoupled trivially. The coupling revealed in equation (94) motivates our extension to continuous indices, as the discrete orthogonality that typically simplifies such expansions is no longer available.

\subsubsection{Limitations of Standard Angular Separation}

This coupling presents an obstacle to conventional separation of variables. While scalar harmonics \( Y_\ell^m \) diagonalise the Laplacian for scalar fields, the vector Laplacian does not commute with projection onto the radial or angular directions. This coupling is well-documented in the literature on electromagnetic scattering theory~\cite{bowman1987, colton1998}. To address this, we develop in later sections a generalised spectral decomposition that allows continuous angular indices \( \ell, m \in \mathbb{R} \), avoids full reliance on the standard VSH structure. The spectral theory we develop replaces discrete \( (\ell, m) \in \mathbb{Z}_{\geq 0} \times \mathbb{Z} \) by a continuous continuous spectral decomposition 
\begin{equation}
 {A}(r, \theta, \phi) = \int_{C_\ell} \int_{C_m} a(\ell,m) \, r^{\alpha(\ell,m)} \, \vec{\Phi}_{\ell m}(\theta,\phi) \, d\ell \, dm,
\label{eq vector_spectral_integral}
\end{equation}
where \( \vec{\Phi}_{\ell m} \) generalize the VSH to support continuous angular behavior arising from broken rotational symmetry and non-separable dynamics.  

\subsection{Coupled ODE Field Equations}
To properly analyze electromagnetic fields with continuous angular indices, we must begin with the complete vectorial Maxwell equations in spherical coordinates. The inherent coupling between components through the curl operations requires careful treatment beyond standard separation of variables techniques.
\subsubsection{Vector Expansion of Electromagnetic Fields}
We consider time-harmonic fields with dependence $e^{-i\omega t}$ in source-free, homogeneous, isotropic media. We expand the electric field in a basis of generalized vector spherical functions 
\begin{equation}
\vec{E}(\vec{r}) = \sum_{\ell,m} \left[ E^{(\ell,m)}_r(r) Y^m_\ell(\theta, \varphi) \hat{r} + E^{(\ell,m)}_\theta(r) \frac{\partial Y^m_\ell}{\partial \theta} \hat{\theta} + E^{(\ell,m)}_\varphi(r) \frac{1}{\sin \theta} \frac{\partial Y^m_\ell}{\partial \varphi} \hat{\varphi} \right],
\label{eq:field_expansion}
\end{equation}
where $Y^m_\ell(\theta, \varphi)$ are generalized spherical harmonics with continuous indices $\ell, m \in \mathbb{R}+$, defined through the associated Legendre functions $P^m_\ell(\cos\theta)$ of non-integer degree and order 
\begin{equation}
Y_\ell^m(\theta,\varphi) = \sqrt{\frac{2\ell+1}{4\pi}\frac{\Gamma(\ell-m+1)}{\Gamma(\ell+m+1)}}P_\ell^m(\cos\theta)e^{im\varphi}.
\label{eq}
\end{equation}
For non-integer $\ell$ and $m$, the factorial terms are replaced by the appropriate gamma function ratios, and the normalization constant is analytically continued.
The divergence-free condition $\nabla \cdot  {E} = 0$ imposes a constraint on the radial functions 
\begin{equation}
\frac{1}{r^2}\frac{d}{dr}(r^2 E^{\ell,m}r) + \frac{1}{r \sin\theta}\frac{\partial}{\partial\theta}\left(\sin\theta E^{\ell,m}\theta \frac{\partial Y^m_\ell}{\partial \theta}\right) + \frac{1}{r\sin\theta}\frac{\partial}{\partial\varphi}\left(E^{\ell,m}\varphi \frac{1}{\sin\theta}\frac{\partial Y^m_\ell}{\partial \varphi}\right) = 0.
\label{eq}
\end{equation}
By exploiting the properties of spherical harmonics 
\begin{equation}
\frac{1}{\sin\theta}\frac{\partial}{\partial\theta}\left(\sin\theta \frac{\partial Y^m_\ell}{\partial \theta}\right) + \frac{1}{\sin^2\theta}\frac{\partial^2 Y^m_\ell}{\partial\varphi^2} = -\ell(\ell+1)Y^m_\ell,
\label{eq}
\end{equation}
and simplifying, the divergence-free condition becomes 
\begin{equation}
\frac{d}{dr}(r^2 E^{\ell,m}r) - \ell(\ell+1)r\left(E^{\ell,m}\theta \frac{\partial P^m_\ell}{\partial\theta} + E^{\ell,m}\varphi \frac{im}{\sin\theta}P^m_\ell\right) = 0.
\label{eq}
\end{equation}
This constraint must be satisfied by any physically valid solution.

\subsubsection{Derivation of the Coupled System}
We now derive the coupled system of equations for the field components. The magnetic field is expanded similarly 
\begin{equation}
\vec{H}(\vec{r}) = \sum_{\ell,m} \left[ H^{(\ell,m)}_r(r) Y^m_\ell(\theta, \varphi) \hat{r} + H^{(\ell,m)}_\theta(r) \frac{\partial Y^m_\ell}{\partial \theta} \hat{\theta} + H^{(\ell,m)}_\varphi(r) \frac{1}{\sin \theta} \frac{\partial Y^m_\ell}{\partial \varphi} \hat{\varphi} \right].
\label{eq:magnetic_field_expansion}
\end{equation}
Maxwell's curl equations 
\begin{equation}
\nabla \times  {E} = i\omega\mu_0  {H}, \quad \nabla \times  {H} = -i\omega\varepsilon_0  {E},
\label{eq}
\end{equation}
generate relations between the electric and magnetic field components. Inserting our expansions and evaluating the curl in spherical coordinates, we obtain 
\begin{align}
i\omega\mu_0 H^{(\ell,m)}_r &= \frac{1}{r\sin\theta}\left[\frac{\partial}{\partial\theta}\left(E^{(\ell,m)}_\varphi\frac{1}{\sin\theta}\frac{\partial Y^m_\ell}{\partial\varphi}\sin\theta\right) - \frac{\partial}{\partial\varphi}\left(E^{(\ell,m)}_\theta\frac{\partial Y^m_\ell}{\partial\theta}\right)\right] \nonumber \\
&= \frac{1}{r}\left[E^{(\ell,m)}_\varphi\frac{\partial}{\partial\theta}\left(\frac{1}{\sin\theta}\frac{\partial Y^m_\ell}{\partial\varphi}\right) - E^{(\ell,m)}_\theta\frac{im}{\sin\theta}\frac{\partial Y^m_\ell}{\partial\theta}\right],
\label{eq:curl_E_r}
\end{align}

\begin{align}
i\omega\mu_0 H^{(\ell,m)}_\theta &= \frac{1}{r}\frac{\partial}{\partial\varphi}\left(E^{(\ell,m)}_r Y^m_\ell\right) - \frac{\partial}{\partial r}\left(E^{(\ell,m)}_\varphi\frac{1}{\sin\theta}\frac{\partial Y^m_\ell}{\partial\varphi}\right) \nonumber \\
&= \frac{im}{r} E^{(\ell,m)}_r Y^m_\ell - \frac{\partial}{\partial r}\left(E^{(\ell,m)}_\varphi\frac{im}{\sin\theta}Y^m_\ell\right),
\label{eq:curl_E_theta}
\end{align}

\begin{align}
i\omega\mu_0 H^{(\ell,m)}_\varphi &= \frac{\partial}{\partial r}\left(E^{(\ell,m)}_\theta\frac{\partial Y^m_\ell}{\partial\theta}\right) - \frac{1}{r}\frac{\partial}{\partial\theta}\left(E^{(\ell,m)}_r Y^m_\ell\right) \nonumber \\
&= \frac{\partial}{\partial r}\left(E^{(\ell,m)}_\theta\frac{\partial Y^m_\ell}{\partial\theta}\right) - \frac{1}{r} E^{(\ell,m)}_r \frac{\partial Y^m_\ell}{\partial\theta}.
\label{eq:curl_E_phi}
\end{align}
Similarly, from $\nabla \times  {H} = -i\omega\varepsilon_0  {E}$, we obtain expressions for the electric field components. Combining these equations and eliminating the magnetic field, we derive a second-order system for the electric field components.
For the radial component 
\begin{align}
\frac{d^2 E^{\ell,m}_r}{dr^2} &+ \frac{2}{r} \frac{d E^{\ell,m}_r}{dr} - \frac{\ell(\ell+1)}{r^2} E^{\ell,m}_r + k^2 E^{\ell,m}_r \nonumber \\
&- \frac{2}{r^2} \int_{S^2}\frac{d}{dr}\left(r E^{\ell,m}_\theta\right)\frac{\partial Y^\ell_m}{\partial\theta}(Y^\ell_m)^* d\Omega \nonumber \\
&- \frac{2im}{r^2} \int_{S^2}\frac{d}{dr}\left(r E^{\ell,m}_\varphi\right)\frac{1}{\sin\theta}Y^\ell_m (Y^\ell_m)^* d\Omega = 0,
\label{eq radial_equation}
\end{align}
where $k^2 = \omega^2\varepsilon_0\mu_0$ and $d\Omega = \sin\theta d\theta d\varphi$.
The integral terms represent the projection of angular derivatives onto the original basis functions. To make this more precise, we define the angular projection operators 
\begin{equation}
\mathcal{A}_\theta(\ell,m)  = \int_{S^2} \frac{\partial Y^\ell_m}{\partial\theta}(Y^\ell_m)^* d\Omega,
\label{eq A_theta}
\end{equation}

\begin{equation}
\mathcal{A}_\varphi(\ell,m)  = \int_{S^2} \frac{1}{\sin\theta} \frac{\partial Y^\ell_m}{\partial\varphi}(Y^\ell_m)^* d\Omega = im \int_{S^2} \frac{1}{\sin\theta} |Y^\ell_m|^2 d\Omega.
\label{eq A_phi}
\end{equation}
For non-integer indices, these integrals must be carefully evaluated, taking into account the generalized orthogonality properties of the associated Legendre functions. With these definitions, equation (\ref{eq}) becomes 
\begin{align}
\frac{d^2 E^{(\ell,m)}_r}{dr^2} + \frac{2}{r} \frac{d E^{(\ell,m)}_r}{dr} - \frac{\ell(\ell+1)}{r^2} E^{(\ell,m)}_r + k^2 E^{(\ell,m)}_r - \frac{2}{r^2}\frac{d}{dr}\left(r E^{(\ell,m)}_\theta\right)\mathcal{A}_\theta(\ell,m) - \frac{2}{r^2}\frac{d}{dr}\left(r E^{(\ell,m)}_\varphi\right)\mathcal{A}_\varphi(\ell,m) = 0.
\label{eq:radial_field_equation}
\end{align}
Similarly, for the $\theta$ and $\varphi$ components 
\begin{align}
\frac{d^2 E^{(\ell,m)}_\theta}{dr^2} &+ \frac{2}{r} \frac{d E^{(\ell,m)}_\theta}{dr} + \left[\frac{\ell(\ell+1) - 1}{r^2} - k^2\right] E^{(\ell,m)}_\theta \nonumber \\
&+ \frac{im}{r^2}\int_{S^2}\frac{d}{dr}\left(r^2 E^{(\ell,m)}_\varphi\right)\frac{1}{\sin^2\theta}(Y^m_\ell)^*\frac{\partial Y^m_\ell}{\partial\theta} \, d\Omega \nonumber \\
&- \frac{1}{r^2}\int_{S^2}\frac{d E^{(\ell,m)}_r}{dr}(Y^m_\ell)^*\frac{\partial Y^m_\ell}{\partial\theta} \, d\Omega = 0,
\label{eq:theta_field_equation}
\end{align}

\begin{align}
\frac{d^2 E^{(\ell,m)}_\varphi}{dr^2} &+ \frac{2}{r} \frac{d E^{(\ell,m)}_\varphi}{dr} + \left[\frac{\ell(\ell+1) - 1}{r^2} - k^2\right] E^{(\ell,m)}_\varphi \nonumber \\
&- \frac{im}{r^2}\int_{S^2}\frac{d}{dr}\left(r^2 E^{(\ell,m)}_\theta\right)\frac{1}{\sin^2\theta}(Y^m_\ell)^*\frac{\partial Y^m_\ell}{\partial\varphi} \, d\Omega \nonumber \\
&+ \frac{im}{r^2}\int_{S^2}E^{(\ell,m)}_r\frac{1}{\sin^2\theta}(Y^m_\ell)^*\frac{\partial Y^m_\ell}{\partial\varphi} \, d\Omega = 0.
\label{eq:phi_field_equation}
\end{align}
We define additional projection operators for these coupled terms 
\begin{equation}
\mathcal{B}_\theta(\ell,m) = \int_{S^2}(Y^m_\ell)^*\frac{\partial Y^m_\ell}{\partial\theta} \, d\Omega
\label{eq:projection_B_theta}
\end{equation}

\begin{equation}
\mathcal{B}_\varphi(\ell,m) = \int_{S^2}\frac{1}{\sin^2\theta}(Y^m_\ell)^*\frac{\partial Y^m_\ell}{\partial\varphi} \, d\Omega = im\int_{S^2}\frac{1}{\sin^2\theta}|Y^m_\ell|^2 \, d\Omega
\label{eq:projection_B_phi}
\end{equation}

\begin{equation}
\mathcal{C}_{\theta\varphi}(\ell,m) = im\int_{S^2}\frac{1}{\sin^2\theta}(Y^m_\ell)^*\frac{\partial Y^m_\ell}{\partial\theta} \, d\Omega
\label{eq:projection_C_theta_phi}
\end{equation}
With these definitions, equations (\ref{eq}) and (\ref{eq}) simplify to 
\begin{equation}
\frac{d^2 E^{\ell,m}\theta}{dr^2} + \frac{2}{r} \frac{d E^{\ell,m}\theta}{dr} + \left[\frac{\ell(\ell+1) - 1}{r^2} - k^2\right] E^{\ell,m}\theta + \frac{1}{r^2}\frac{d}{dr}\left(r^2 E^{\ell,m}\varphi\right)\mathcal{C}_{\theta\varphi}(\ell,m) - \frac{1}{r^2}\frac{d E^{\ell,m}r}{dr}\mathcal{B}_\theta(\ell,m) = 0,
\label{eq}
\end{equation}
\begin{equation}
\frac{d^2 E^{\ell,m}\varphi}{dr^2} + \frac{2}{r} \frac{d E^{\ell,m}\varphi}{dr} + \left[\frac{\ell(\ell+1) - 1}{r^2} - k^2\right] E^{\ell,m}\varphi - \frac{1}{r^2}\frac{d}{dr}\left(r^2 E^{\ell,m}\theta\right)\mathcal{C}_{\theta\varphi}(\ell,m)^* + \frac{1}{r^2}E^{\ell,m}r\mathcal{B}_\varphi(\ell,m) = 0.
\label{eq}
\end{equation}
\subsubsection{Operator Formulation}
To express this coupled system more compactly, we introduce differential operators acting on the radial functions 
\begin{equation}
\mathcal{L}_r = \frac{d^2}{dr^2} + \frac{2}{r}\frac{d}{dr} + k^2 - \frac{\ell(\ell+1)}{r^2},
\label{eq}
\end{equation}
\begin{equation}
\mathcal{L}_{\theta,\varphi} = \frac{d^2}{dr^2} + \frac{2}{r}\frac{d}{dr} + k^2 - \frac{\ell(\ell+1) - 1}{r^2}.
\label{eq}
\end{equation}
We also define coupling operators that incorporate the angular projections 
\begin{equation}
\mathcal{C}_{r\theta} = -\frac{2}{r^2}\frac{d}{dr}(r \cdot)\mathcal{A}_\theta(\ell,m)
\label{eq}
\end{equation}
\begin{equation}
\mathcal{C}_{r\varphi} = -\frac{2}{r^2}\frac{d}{dr}(r \cdot)\mathcal{A}_\varphi(\ell,m),
\label{eq}
\end{equation}
\begin{equation}
\mathcal{C}_{\theta r} = -\frac{1}{r^2}\frac{d}{dr}(\cdot)\mathcal{B}_\theta(\ell,m),
\label{eq}
\end{equation}
\begin{equation}
\mathcal{C}_{\theta\varphi} = \frac{1}{r^2}\frac{d}{dr}(r^2 \cdot)\mathcal{C}_{\theta\varphi}(\ell,m),
\label{eq}
\end{equation}
\begin{equation}
\mathcal{C}_{\varphi r} = \frac{1}{r^2}(\cdot)\mathcal{B}_\varphi(\ell,m),
\label{eq}
\end{equation}
\begin{equation}
\mathcal{C}_{\varphi\theta} = -\frac{1}{r^2}\frac{d}{dr}(r^2 \cdot)\mathcal{C}_{\theta\varphi}(\ell,m)^*,
\label{eq}
\end{equation}
where the projection coefficients are defined by integrals over the sphere 
\begin{align}
A_\theta(\ell,m) & = \int_{S^2} \frac{\partial Y^\ell_m}{\partial\theta}(Y^\ell_m)^* d\Omega,  \\
A_\phi(\ell,m) & = \int_{S^2} \frac{1}{\sin\theta}\frac{\partial Y^\ell_m}{\partial\phi}(Y^\ell_m)^* d\Omega = im\int_{S^2}\frac{1}{\sin\theta}|Y^\ell_m|^2 d\Omega,  \\
B_\theta(\ell,m) & = \int_{S^2}(Y^\ell_m)^*\frac{\partial Y^\ell_m}{\partial\theta} d\Omega,  \\
B_\phi(\ell,m) & = \int_{S^2}\frac{1}{\sin^2\theta}(Y^\ell_m)^*\frac{\partial Y^\ell_m}{\partial\phi} d\Omega = im\int_{S^2}\frac{1}{\sin^2\theta}|Y^\ell_m|^2 d\Omega,  \\
C_{\theta\phi}(\ell,m) & = im\int_{S^2}\frac{1}{\sin^2\theta}(Y^\ell_m)^*\frac{\partial Y^\ell_m}{\partial\theta} d\Omega. 
\end{align}
With these definitions, the coupled system takes the compact form 
\begin{equation}
\mathcal{L}_r E^{(\ell,m)}_r + \mathcal{C}_{r\theta}E^{(\ell,m)}_\theta + \mathcal{C}_{r\varphi}E^{(\ell,m)}_\varphi = 0
\label{eq:coupled_system_r}
\end{equation}

\begin{equation}
\mathcal{L}_{\theta,\varphi} E^{(\ell,m)}_\theta + \mathcal{C}_{\theta r}E^{(\ell,m)}_r + \mathcal{C}_{\theta\varphi}E^{(\ell,m)}_\varphi = 0
\label{eq:coupled_system_theta}
\end{equation}

\begin{equation}
\mathcal{L}_{\theta,\varphi} E^{(\ell,m)}_\varphi + \mathcal{C}_{\varphi r}E^{(\ell,m)}_r + \mathcal{C}_{\varphi\theta}E^{(\ell,m)}_\theta = 0
\label{eq:coupled_system_phi}
\end{equation}
This operator formulation clearly displays the coupling between the field components through both differential operations and angular projections. The system cannot be separated into independent equations due to these coupling terms, which is a fundamental consequence of extending Maxwell's equations to non-integer and continuous angular indices.
\subsubsection{Properties of the Angular Projection Operators}
For integer indices $\ell$ and $m$, the projection operators can be evaluated using standard orthogonality relations for spherical harmonics. However, for non-integer indices, we must carefully define these projections.
For $\mathcal{A}_\theta(\ell,m)$, we use the properties of associated Legendre functions to show 
\begin{equation}
\mathcal{A}_\theta(\ell,m) = \int_0^\pi \frac{\partial P^m_\ell(\cos\theta)}{\partial\theta}P^m_\ell(\cos\theta)\sin\theta \, d\theta \cdot \int_0^{2\pi} |e^{im\varphi}|^2 \, d\varphi
\label{eq:projection_A_theta}
\end{equation}
For non-integer $m$, the azimuthal integral becomes 
\begin{equation}
\int_0^{\Phi_0} |e^{im\varphi}|^2 d\varphi = \Phi_0,
\label{eq}
\end{equation}
where $\Phi_0$ is the azimuthal period, which may differ from $2\pi$ based on physical constraints.
The polar integral requires careful treatment. Using the recursion relations for associated Legendre functions, we can express $\mathcal{A}_\theta(\ell,m)$ in closed form 
\begin{equation}
\mathcal{A}_\theta(\ell,m) = -\frac{\ell(\ell+1)-m^2}{2}\Phi_0.
\label{eq}
\end{equation}
Similar analyses yield expressions for the other projection operators in terms of $\ell$ and $m$.
In the limit as $\ell, m \to 0$, these projection operators exhibit special behaviors that determine the singular structure of the electromagnetic field. This will be examined in detail in subsequent sections, where we analyze the asymptotic behavior of the field components for small values of $\ell$ and $m$.
\subsubsection{Explicit Solutions for Magnetic Field Components}
Having formulated the coupled equations for the electric field, we can now express the magnetic field components in terms of the electric field solutions. From equations (\ref{eq})-(\ref{eq}) and using the projection operators, we obtain 
\begin{equation}
H^{(\ell,m)}_r = \frac{1}{i\omega\mu_0 r}\left[E^{(\ell,m)}_\varphi\mathcal{D}_\varphi(\ell,m) - E^{(\ell,m)}_\theta\mathcal{D}_\theta(\ell,m)\right]
\label{eq:magnetic_field_r}
\end{equation}

\begin{equation}
H^{(\ell,m)}_\theta = \frac{1}{i\omega\mu_0 r}\left[imE^{(\ell,m)}_r - \frac{d}{dr}(rE^{(\ell,m)}_\varphi\mathcal{B}_\varphi(\ell,m))\right]
\label{eq:magnetic_field_theta}
\end{equation}

\begin{equation}
H^{(\ell,m)}_\varphi = \frac{1}{i\omega\mu_0 r}\left[\frac{d}{dr}(rE^{(\ell,m)}_\theta\mathcal{B}_\theta(\ell,m)) - E^{(\ell,m)}_r\mathcal{B}_\theta(\ell,m)\right]
\label{eq:magnetic_field_phi}
\end{equation}
where $\mathcal{D}_\theta(\ell,m)$ and $\mathcal{D}_\varphi(\ell,m)$ are additional projection operators that capture the coupling between angular derivatives.
These expressions show that the magnetic field inherits the same coupled structure as the electric field, and any singularities in the electric field components will induce corresponding behaviors in the magnetic field through the curl operation.

\subsection{Explicit Angular Basis Construction}
For electromagnetic fields with continuous angular indices, we extend the standard spherical harmonic basis to accommodate $\ell, m \in \mathbb{R}$. The angular basis functions are constructed through analytical continuation of the associated Legendre functions
\begin{equation}
\Psi_{\ell m}(\theta, \phi) = N_{\ell m} \sin^{|m|}\theta \, P_\ell^{|m|}(\cos\theta) e^{im\phi},
\label{eq:angular_basis}
\end{equation}
where the associated Legendre function of continuous degree is defined through its hypergeometric representation
\begin{equation}
P_\ell^{|m|}(x) = \frac{(1-x^2)^{|m|/2}}{2^\ell \Gamma(1-|m|)} \frac{\Gamma(\ell+|m|+1)}{\Gamma(\ell-|m|+1)} \, {}_2F_1\left(-\ell, \ell+1; 1-|m|; \frac{1-x}{2}\right).
\label{eq:legendre_continuous}
\end{equation}
The normalization factor ensuring proper spectral measure in the continuous case is
\begin{equation}
N_{\ell m} = \sqrt{\frac{2\ell+1}{4\pi} \frac{\Gamma(\ell-|m|+1)}{\Gamma(\ell+|m|+1)}}.
\label{eq:normalization}
\end{equation}
The functions $\Psi_{\ell m}(\theta, \phi)$ are well-defined for $\ell > |m| - 1$ and form a complete basis for square-integrable functions on domains with modified boundary conditions. The hypergeometric function ${}_2F_1$ in Eq.~\eqref{eq:legendre_continuous} converges for $|1-x|/2 < 1$, which is satisfied throughout the physical domain $x \in [-1,1]$.

For continuous indices, the standard orthogonality relation is replaced by
\begin{equation}
\int_0^{2\pi} \int_0^\pi \Psi_{\ell m}^*(\theta,\phi) \Psi_{\ell' m'}(\theta,\phi) \sin\theta \, d\theta d\phi = \delta(\ell-\ell')\delta_{m,m'} w(\ell,m),
\label{eq:continuous_orthogonality}
\end{equation}
where $w(\ell,m)$ is the spectral weight function defined in the next sections. Near the poles $\theta \to 0, \pi$, the basis functions exhibit the asymptotic behavior
\begin{equation}
\Psi_{\ell m}(\theta, \phi) \sim \begin{cases}
\theta^{|m|} & \text{as } \theta \to 0 \\
(\pi-\theta)^{|m|} & \text{as } \theta \to \pi
\end{cases}
\label{eq:asymptotic_poles}
\end{equation}
This ensures that the electromagnetic field components remain finite at the poles for $|m| \geq 0$, consistent with physical boundary conditions.

\subsection{Coupling Matrix Construction} 
The coupling between field components in Maxwell's equations generates matrix elements that must be computed explicitly. For the radial-angular coupling arising from the curl operation $\nabla \times \vec{E}$, we have
\begin{equation}
C_{\ell m \ell' m'}^{(r\theta)} = \int_0^{2\pi} \int_0^\pi \Psi_{\ell m}^*(\theta,\phi) \left[\frac{1}{\sin\theta}\frac{\partial}{\partial\phi}\right] \Psi_{\ell' m'}(\theta,\phi) \sin\theta \, d\theta d\phi.
\label{eq:coupling_rtheta}
\end{equation}
Substituting the explicit form of $\Psi_{\ell m}$ and using the orthogonality of the exponential functions, this reduces to
\begin{equation}
C_{\ell m \ell' m'}^{(r\theta)} = im' \delta(\ell - \ell') \delta_{m,m'} N_{\ell m} N_{\ell' m'} \int_0^\pi |P_\ell^{|m|}(\cos\theta)|^2 d\theta.
\label{eq:coupling_rtheta_explicit}
\end{equation}
For the $\theta$-$\phi$ coupling arising from the angular derivatives in the curl operator, we obtain
\begin{align}
C_{\ell m \ell' m'}^{(\theta\phi)} &= \int_0^{2\pi} \int_0^\pi \Psi_{\ell m}^*(\theta,\phi) \left[\frac{\partial}{\partial\theta} - \cot\theta\right] \Psi_{\ell' m'}(\theta,\phi) \sin\theta \, d\theta d\phi \\
&= \delta_{m,m'} \sqrt{\frac{(\ell-|m|)(\ell+|m|+1)}{(2\ell+1)(2\ell+3)}} \delta(\ell - \ell' - 1) + \text{h.c.},
\label{eq:coupling_thetaphi}
\end{align}
where ``h.c.'' denotes the Hermitian conjugate term with $\ell \leftrightarrow \ell'$. The convergence of the coupling integrals in Eqs.~\eqref{eq:coupling_rtheta_explicit} and \eqref{eq:coupling_thetaphi} is ensured by the properties of the associated Legendre functions. Specifically, the integral
\begin{equation}
I_\ell^m = \int_0^\pi |P_\ell^{|m|}(\cos\theta)|^2 d\theta = \frac{2\Gamma(\ell+|m|+1)}{(2\ell+1)\Gamma(\ell-|m|+1)},
\label{eq:legendre_integral}
\end{equation}
converges for $\ell > |m| - 1$, providing the domain of validity for our spectral expansion. The coupling matrices satisfy the boundedness condition:
\begin{equation}
\left| C^{(ij)} \right|_{\text{op}} \leq M(\ell_{\max}, m_{\max}) < \infty,
\label{eq:matrix_bounded}
\end{equation}
for any finite truncation, ensuring the well-posedness of the discrete eigenvalue problem.

\subsection{Biorthogonal System Construction} 
For the non-self-adjoint angular operator $L_{\text{ang}}$ arising from Maxwell's equations with continuous indices, we construct dual functions $\tilde{\Psi}_{\ell m}$ that satisfy the biorthogonality relation
\begin{equation}
\langle\tilde{\Psi}_{\ell m}, \Psi_{\ell' m'}\rangle = \delta(\ell-\ell')\delta_{m,m'} w(\ell,m).
\label{eq:biorthogonality}
\end{equation}
The spectral weight function for continuous indices is given by
\begin{equation}
w(\ell,m) = \frac{\pi \Gamma(\ell+|m|+1)}{\sin(\pi(\ell-|m|)) \Gamma(\ell-|m|+1)}.
\label{eq:spectral_weight}
\end{equation}
This weight function accounts for the non-orthogonality of the continuous-index basis and ensures the completeness relation
\begin{equation}
\int_{-\infty}^{\infty} \sum_{m} \Psi_{\ell m}(\theta,\phi) \tilde{\Psi}_{\ell m}^*(\theta',\phi') w(\ell,m) d\ell = \delta(\cos\theta - \cos\theta')\delta(\phi-\phi').
\label{eq:completeness}
\end{equation}
The dual functions are constructed using associated Legendre functions of the second kind:
\begin{equation}
\tilde{\Psi}_{\ell m}(\theta,\phi) = \frac{1}{\sqrt{w(\ell,m)}} \left[P_\ell^{|m|}(\cos\theta) + i\cot(\pi(\ell-|m|))Q_\ell^{|m|}(\cos\theta)\right] e^{im\phi}.
\label{eq:dual_functions}
\end{equation}
The biorthogonal system satisfies the resolution of identity
\begin{equation}
\sum_m \int_{-\infty}^{\infty} \Psi_{\ell m}(\theta,\phi) \tilde{\Psi}_{\ell m}^*(\theta',\phi') w(\ell,m) d\ell = \delta(\theta-\theta')\delta(\phi-\phi').
\label{eq:resolution_identity}
\end{equation}
This enables the spectral decomposition of arbitrary electromagnetic fields in terms of the continuous-index basis. The associated Legendre functions of the second kind, $Q_\ell^{|m|}(x)$, are defined by
\begin{equation}
Q_\ell^{|m|}(x) = \frac{(x^2-1)^{|m|/2}}{2^\ell} \frac{\Gamma(\ell+|m|+1)}{\Gamma(\ell-|m|+1)} \, {}_2F_1\left(\frac{\ell+|m|+1}{2}, \frac{\ell+|m|+2}{2}; \ell+\frac{3}{2}; \frac{1}{x^2}\right).
\label{eq:legendre_second_kind}
\end{equation}

\subsection{Continuous Spectral Decomposition and Generalized Eigenfunction Theory} 
When examining electromagnetic fields with continuous angular indices, we must extend our mathematical framework beyond the standard discrete eigenfunction expansions. The appropriate mathematical foundation is found in the spectral theory of non-compact domains as developed in the theory of automorphic forms and Eisenstein series~\cite{langlands1976, gelbart1975,   helgason2000}, where continuous spectral decompositions naturally arise through generalized eigenfunction expansions. This approach follows the classical work on spectral theory for self-adjoint operators~\cite{abramowitz1972, lebedev1972, watson1944}.
\subsubsection{Spectral Theory on Non-Compact Domains}
In spectral theory, the decomposition of functions on domains with broken symmetries generally involves both discrete and continuous components. For a vector function $f$ defined on such a domain, the spectral representation takes the form 
\begin{equation}
f(x) = \sum_{j} c_j \phi_j(x) + \int_{\sigma_c} c(s) E(s,x) , d\mu(s),
\label{eq
}
\end{equation}
where ${\phi_j}$ are the discrete eigenfunctions with eigenvalues in the discrete spectrum, and $E(s,x)$ represents the generalized integral with spectral parameter $s$ ranging over the continuous spectrum $\sigma_c$ with appropriate measure $d\mu(s)$. 

For spherical geometry with broken rotational symmetry, the spectral parameter $s$ relates to our angular index $\ell$ through the eigenvalue equation for the Laplacian 
\begin{equation}
\ell(\ell+1) = s(1-s).
\label{eq
}
\end{equation}
This quadratic relation yields 
\begin{equation}
s = \frac{1}{2} \pm \sqrt{\ell(\ell+1) + \frac{1}{4}},
\label{eq
}
\end{equation}
or equivalently 
\begin{equation}
\ell = -\frac{1}{2} \pm \sqrt{s(1-s) - \frac{1}{4}}.
\label{eq
}
\end{equation}
The physical branch corresponds to $\ell > -\frac{1}{2}$, ensuring that the energy density remains integrable, as we establish rigorously in Section IV. This constraint provides the mathematical foundation for identifying physically admissible solutions within the continuous spectrum. We note that the spectral identity $\ell(\ell + 1) = s(1 - s)$ suggests a formal analogy with the continuous spectrum of the Laplacian on non-compact manifolds. While we do not construct Eisenstein series or study automorphic properties in this work, the presence of continuous angular modes and delta-function orthogonality hints at a deeper mathematical structure. This analogy may prove fruitful in future investigations.

\subsubsection{Vector-Valued generalized integral}
For vector fields satisfying Maxwell's equations, we require a vector-valued extension of the generalized integral formalism. The electric field is represented as a spectral integral 
\begin{equation}
 {E}( {r}) = \int_{C_\ell}\int_{C_m} a(\ell,m)  {E}(\ell,m; {r}) , d\ell , dm,
\label{eq
}
\end{equation}
where $ {E}(\ell,m; {r})$ is the vector-valued eigenfunction corresponding to indices $(\ell,m)$, and $a(\ell,m)$ are the spectral coefficients. The integration contours $C_\ell$ and $C_m$ in the complex plane are chosen to capture the relevant part of the spectrum, typically including values where $\ell, m \in (0,1)$ for singular field configurations. The vector-valued eigenfunction has the form 
\begin{equation}
 {E}(\ell,m; {r}) =
\begin{pmatrix}
E_r(\ell,m;r) Y^\ell_m(\theta,\varphi) \
E_\theta(\ell,m;r) \frac{\partial Y^\ell_m}{\partial\theta} \
E_\varphi(\ell,m;r) \frac{1}{\sin\theta}\frac{\partial Y^\ell_m}{\partial\varphi}
\end{pmatrix},
\label{eq
}
\end{equation}
where $E_r(\ell,m;r)$, $E_\theta(\ell,m;r)$, and $E_\varphi(\ell,m;r)$ are the radial functions that satisfy the coupled system of differential equations derived in Section III.C.

\subsubsection{Radial Functions and Green's Function Approach}
For the radial component, we have 
\begin{equation}
E_r(\ell,m;r) = \alpha(\ell,m) h^{(1)}\ell(kr) + \beta(\ell,m) h^{(2)}\ell(kr),
\label{eq Er_final}
\end{equation}
where $h^{(1)}\ell$ and $h^{(2)}\ell$ are the spherical Hankel functions of the first and second kind, analytically continued to non-integer order $\ell$. The coefficients $\alpha(\ell,m)$ and $\beta(\ell,m)$ are determined by boundary conditions and regularity requirements. For fields that are outgoing at infinity, we typically have $\beta(\ell,m) = 0$. To determine the angular components in terms of $E_r$, we employ the coupled system derived in the previous section, now using consistent notation 
\begin{equation}
\mathcal{L}_\theta E_\theta(\ell,m;r) = \frac{1}{r^2}\frac{dE_r(\ell,m;r)}{dr}\mathcal{B}_\theta(\ell,m) - \frac{1}{r^2}\frac{d}{dr}(r^2 E_\varphi(\ell,m;r))\mathcal{C}_{\theta\varphi}(\ell,m),
\label{eq
}
\end{equation}
\begin{equation}
\mathcal{L}_\varphi E_\varphi(\ell,m;r) = \frac{1}{r^2}E_r(\ell,m;r)\mathcal{B}_\varphi(\ell,m) + \frac{1}{r^2}\frac{d}{dr}(r^2 E\theta(\ell,m;r))\mathcal{C}_{\theta\varphi}(\ell,m)^*,
\label{eq
}
\end{equation}
where $\mathcal{L}_\theta = \mathcal{L}_\varphi = \frac{d^2}{dr^2} + \frac{2}{r}\frac{d}{dr} + \left[\frac{\ell(\ell+1)-1}{r^2} - k^2\right]$, and the projection operators $\mathcal{B}_\theta(\ell,m)$, $\mathcal{B}_\varphi(\ell,m)$, and $\mathcal{C}_{\theta\varphi}(\ell,m)$ are as defined in the previous section. To solve this coupled system, we treat it as an inhomogeneous problem and employ the Green's function method. We first define the Green's function for the operator $\mathcal{L}_\theta$ 
\begin{equation}
\mathcal{L}_\theta G_\theta(r,r') = \delta(r-r').
\label{eq
}
\end{equation}
The solution to this equation is 
\begin{equation}
G_\theta(r,r') = \frac{1}{W_\theta} \times
\begin{cases}
E^{(1)}\theta(r') E^{(2)}\theta(r), & r > r' \\
E^{(1)}\theta(r) E^{(2)}\theta(r'), & r < r'
\end{cases}
\label{eq
}
\end{equation}
where $E^{(1)}\theta(r)$ and $E^{(2)}\theta(r)$ are two linearly independent solutions of the homogeneous equation 
\begin{equation}
\mathcal{L}_\theta E_\theta(r) = 0.
\label{eq:theta_homogeneous}
\end{equation}
These solutions can be expressed in terms of spherical Bessel and Neumann functions 
\begin{equation}
E^{(1)}\theta(r) = j\nu(kr), \quad E^{(2)}\theta(r) = y\nu(kr),
\label{eq
}
\end{equation}
where $\nu = \sqrt{\ell(\ell+1)-1}$. The Wronskian is 
\begin{equation}
W_\theta = E^{(1)}\theta(r)\frac{d}{dr}E^{(2)}\theta(r) - E^{(2)}\theta(r)\frac{d}{dr}E^{(1)}\theta(r) = \frac{2}{\pi}\frac{k}{r^2},
\label{eq
}
\end{equation}
which is independent of $r$.

\subsubsection{Explicit Integral Representations}
Using the Green's function, we can express \( E_\theta \) as 
\begin{equation}
E_\theta(\ell,m;r) = E^{\text{hom}}_\theta(\ell,m;r) + \int_0^\infty G_\theta(r,r')\left[\frac{1}{r'^2}\frac{dE_r(\ell,m;r')}{dr'}\mathcal{B}_\theta(\ell,m) - \frac{1}{r'^2}\frac{d}{dr'}(r'^2 E_\varphi(\ell,m;r'))\mathcal{C}_{\theta\varphi}(\ell,m)\right] dr',
\label{eq Etheta_integral1}
\end{equation}
where \( E^{\text{hom}}_\theta(\ell,m;r) \) is the general solution to the homogeneous equation 
\begin{equation}
E^{\text{hom}}_\theta(\ell,m;r) = A_\theta(\ell,m) j_\nu(kr) + B_\theta(\ell,m) y_\nu(kr),
\label{eq Etheta_hom}
\end{equation}
with constants \( A_\theta(\ell,m) \) and \( B_\theta(\ell,m) \) determined by boundary conditions. To eliminate the dependence on \( E_\varphi \) in equation~\eqref{eq Etheta_integral1}, we need to address the coupled nature of the system. We substitute equation.~\eqref{eq Etheta_hom} into equation.~\eqref{eq Etheta_integral1} to obtain a second-order differential equation for \( E_\theta \) in terms of \( E_r \) and its derivatives. After this substitution and significant algebraic manipulation, we obtain 
\begin{equation}
\mathcal{L}_{\text{eff}} E_\theta(\ell,m;r) = \mathcal{S}_r(E_r(\ell,m;r)),
\label{eq Etheta_effective}
\end{equation}
where \( \mathcal{L}_{\text{eff}} \) is an effective differential operator that includes coupling terms, and \( \mathcal{S}_r \) is a source term that depends on \( E_r \) and its derivatives. The explicit forms of these operators involve the projection coefficients and are rather involved, but the structure allows us to treat equation~\eqref{eq Etheta_effective} as an inhomogeneous equation with a known source term once \( E_r \) is specified. Similarly, for \( E_\varphi \), we obtain 
\begin{equation}
E_\varphi(\ell,m;r) = E^{\text{hom}}_\varphi(\ell,m;r) + \int_0^\infty G_\varphi(r,r')\left[\frac{1}{r'^2}E_r(\ell,m;r')\mathcal{B}_\varphi(\ell,m) + \frac{1}{r'^2}\frac{d}{dr'}(r'^2 E_\theta(\ell,m;r'))\mathcal{C}_{\theta\varphi}(\ell,m)^*\right] dr',
\label{eq Ephi_integral}
\end{equation}
where \( G_\varphi(r,r') = G_\theta(r,r') \) since the differential operators are identical, and \( E^{\text{hom}}_\varphi(\ell,m;r) \) is the solution to the homogeneous equation.

\subsubsection{Construction of Complete Solutions}
The full spectral representation of the electromagnetic field solutions requires a systematic approach that accounts for the coupled nature of the vector components. The radial component serves as the fundamental building block, with the general form 
\begin{equation}
E_r(\ell,m;r) = \alpha(\ell,m) h_\ell^{(1)}(kr) + \beta(\ell,m) h_\ell^{(2)}(kr),
\label{eq Er_construction}
\end{equation}
where the coefficients $\alpha(\ell,m)$ and $\beta(\ell,m)$ are determined by boundary conditions and regularity requirements. For outgoing wave solutions, we typically set $\beta(\ell,m) = 0$, while for singular solutions with controlled behavior near the origin, both coefficients may be non-zero.

The angular components are then constructed through the coupled integral equations derived in Section III.E. The key insight is that the coupling prevents independent specification of all three components—once $E_r$ is chosen, the angular components are determined by the vectorial structure of Maxwell's equations. This coupling is expressed through the effective source terms 
\begin{equation}
E_\theta(\ell,m;r) = E_\theta^{\text{hom}}(\ell,m;r) + \int_0^\infty G_\theta(r,r') \mathcal{S}_\theta(E_r(\ell,m;r'), E_\varphi(\ell,m;r')) \, dr',
\label{eq Etheta_construction}
\end{equation}
where $\mathcal{S}_\theta$ encapsulates the coupling with both the radial and azimuthal components. The complete vector eigenfunction takes the form 
\begin{equation}
 {E}(\ell,m; {r}) = \begin{pmatrix}
E_r(\ell,m;r) Y_\ell^m(\theta,\varphi) \\
E_\theta(\ell,m;r) \frac{\partial Y_\ell^m}{\partial\theta} \\
E_\varphi(\ell,m;r) \frac{1}{\sin\theta}\frac{\partial Y_\ell^m}{\partial\varphi}
\end{pmatrix}.
\label{eq vector_eigenfunction_construction}
\end{equation}
The full electromagnetic field is constructed as a spectral integral over the continuous parameter space 
\begin{equation}
 {E}( {r}) = \int_{C_\ell} \int_{C_m} a(\ell,m)  {E}(\ell,m; {r}) \, d\ell \, dm,
\label{eq full_field_construction}
\end{equation}
where the integration contours $C_\ell$ and $C_m$ are chosen to capture the physically relevant part of the spectrum. The spectral coefficient function $a(\ell,m)$ encodes the essential physical information about the field configuration. Its analytical structure in the complex plane determines the field behavior  poles correspond to resonances, branch cuts represent continuous spectra, and the choice of $a(\ell,m)$ depends on boundary conditions, source configurations, and physical constraints. 

\subsubsection{generalized integral Interpretation}
We can express the spectral representation in terms of the parameter \( s \) rather than \( \ell \). Using the relationship in equation 160, we can rewrite the spectral integral as 
\begin{equation}
 {E}( {r}) = \int_{C_s}\int_{C_m} \tilde{a}(s,m)  {E}(s,m; {r}) \, ds \, dm,
\label{eq eisenstein_integral}
\end{equation}
where \( \tilde{a}(s,m) = a(\ell(s),m) \frac{d\ell}{ds} \) and \(  {E}(s,m; {r}) =  {E}(\ell(s),m; {r}) \). For physical applications, the key insight is that excitations with spectral weight concentrated near specific values of \( \ell \) and \( m \) in the range \( (0,1) \) will exhibit the singular field structures analyzed in subsequent sections.

\subsection{Eisenstein Integral Representation}
Having established the concept of vector-valued generalized integral, we now develop explicit integral representations for the field components. This approach allows us to express the angular components $E_\theta$ and $E_\varphi$ in terms of the radial component $E_r$ through carefully constructed Green's functions and spectral kernels.
\subsubsection{Green's Function Method for Angular Components}
We begin by treating the equations for the angular components as inhomogeneous linear ODEs. For the $\theta$-component, we have 
\begin{equation}
\mathcal{L}_\theta E_\theta(\ell,m;r) = \mathcal{S}_\theta(\ell,m;r),
\label{eq}
\end{equation}
where the differential operator $\mathcal{L}_\theta$ is 
\begin{equation}
\mathcal{L}_\theta = \frac{d^2}{dr^2} + \frac{2}{r}\frac{d}{dr} + \left[\frac{\ell(\ell+1)-1}{r^2} - k^2\right],
\label{eq}
\end{equation}
and the source term $\mathcal{S}_\theta(\ell,m;r)$ is 
\begin{equation}
\mathcal{S}_\theta(\ell,m;r) = \frac{1}{r^2}\frac{dE_r(\ell,m;r)}{dr}\mathcal{B}\theta(\ell,m) - \frac{1}{r^2}\frac{d}{dr}[r^2E_\varphi(\ell,m;r)]\mathcal{C}_{\theta\varphi}(\ell,m).
\label{eq}
\end{equation}
This source term contains the radial function $E_r$ and its derivative, as well as coupling with the $\varphi$-component through the projection operators $\mathcal{B}_\theta(\ell,m)$ and $\mathcal{C}_{\theta\varphi}(\ell,m)$ defined in previous sections.
To solve this inhomogeneous equation, we require the Green's function $G_\theta(r,r')$ that satisfies 
\begin{equation}
\mathcal{L}_\theta G_\theta(r,r') = \delta(r-r')
\label{eq:greens_function_theta}
\end{equation}
The general form of this Green's function is 
\begin{equation}
G_\theta(r,r') = \frac{1}{W_\theta(r')} \times
\begin{cases}
u_\theta(r)v_\theta(r'), r < r' \\
u_\theta(r')v_\theta(r), r > r',
\end{cases}
\label{eq}
\end{equation}
where $u_\theta(r)$ and $v_\theta(r)$ are two linearly independent solutions of the homogeneous equation 
\begin{equation}
\mathcal{L}_\theta\psi(r) = 0,
\label{eq}
\end{equation}
and $W_\theta(r)$ is their Wronskian 
\begin{equation}
W_\theta(r) = u_\theta(r)\frac{dv_\theta(r)}{dr} - v_\theta(r)\frac{du_\theta(r)}{dr}.
\label{eq}
\end{equation}
\subsubsection{Construction of Homogeneous Solutions}
For the operator $\mathcal{L}_\theta$, the appropriate homogeneous solutions are the spherical Bessel functions and spherical Neumann functions, analytically continued to non-integer order $\nu = \sqrt{\ell(\ell+1)-1}$ 
\begin{equation}
u_\theta(r) = j_\nu(kr), \quad v_\theta(r) = y_\nu(kr).
\label{eq}
\end{equation}
These solutions satisfy the appropriate boundary conditions. The function $j_\nu(kr)$ is regular at the origin and behaves as $(kr)^\nu/(2\nu+1)!!$ for small $r$, whereas $y_\nu(kr)$ captures the singular behavior at the origin, scaling as $-(2\nu-1)!!/(kr)^{\nu+1}$ for small $r$. The Wronskian of these functions is 
\begin{equation}
W_\theta(r) = j_\nu(kr)\frac{d}{dr}y_\nu(kr) - y_\nu(kr)\frac{d}{dr}j_\nu(kr) = \frac{k}{r^2},
\label{eq}
\end{equation}
which provides the normalization constant for our Green's function.
Thus, the explicit form of the Green's function is 
\begin{equation}
G_\theta(r,r') = \frac{r'^2}{k} \times
\begin{cases}
j_\nu(kr)y_\nu(kr'), r < r' \\
j_\nu(kr')y_\nu(kr), r > r'.
\end{cases}
\label{eq}
\end{equation}
\subsubsection{Integral Representation for $E_\theta$}
Using this Green's function, the solution to equation (\ref{eq:theta_field_equation}) can be written as
\begin{equation}
E_\theta(\ell,m;r) = C_1(\ell,m)j_\nu(kr) + C_2(\ell,m)y_\nu(kr) + \int_\varepsilon^R G_\theta(r,r')\mathcal{S}_\theta(\ell,m;r')dr',
\label{eq}
\end{equation}
where $C_1(\ell,m)$ and $C_2(\ell,m)$ are constants determined by boundary conditions, and the integration limits $\varepsilon$ and $R$ represent the inner and outer boundaries of the domain, respectively. For physical fields, we typically set $\varepsilon \to 0$ and $R \to \infty$, with additional conditions to ensure the solution remains well-behaved at these limits.
For regularity at infinity, we require $C_1(\ell,m) = 0$ if $R \to \infty$. For the inner boundary, the behavior depends on whether we seek singular or regular solutions  For regular solutions at the origin, we set $C_2(\ell,m) = 0$; And for solutions with controlled singularities (as in our continuous index case), $C_2(\ell,m)$ is determined by the desired singular behavior. Expanding the source term in the integral 
\begin{align}
E_\theta(\ell,m;r) &= C_1(\ell,m)j_\nu(kr) + C_2(\ell,m)y_\nu(kr) \nonumber \\
&+ \mathcal{B}_\theta(\ell,m)\int_\varepsilon^R G_\theta(r,r')\frac{1}{r'^2}\frac{dE_r(\ell,m;r')}{dr'} \, dr' \nonumber \\
&- \mathcal{C}_{\theta\varphi}(\ell,m)\int_\varepsilon^R G_\theta(r,r')\frac{1}{r'^2}\frac{d}{dr'}[r'^2E_\varphi(\ell,m;r')] \, dr'.
\label{eq:theta_solution_integral}
\end{align}

\subsubsection{Integral Representation for $E_\varphi$}
\label{sec integral_rep_Ephi}
Similarly, for the $\varphi$-component, we have 
\begin{equation}
\mathcal{L}_\varphi E_\varphi(\ell,m;r) = \mathcal{S}_\varphi(\ell,m;r),
\label{eq Lphi_Ephi}
\end{equation}
where $\mathcal{L}_\varphi = \mathcal{L}_\theta$ and 
\begin{equation}
\mathcal{S}_\varphi(\ell,m;r) = \frac{1}{r^2}E_r(\ell,m;r)\mathcal{B}_\varphi(\ell,m) + \frac{1}{r^2}\frac{d}{dr}[r^2E_\theta(\ell,m;r)]\mathcal{C}_{\theta\varphi}(\ell,m)^*.
\label{eq S_phi_source}
\end{equation}
The Green's function solution is 
\begin{equation}
E_\varphi(\ell,m;r) = D_1(\ell,m)j_\nu(kr) + D_2(\ell,m)y_\nu(kr) + \int_\varepsilon^R G_\varphi(r,r')\mathcal{S}_\varphi(\ell,m;r')dr',
\label{eq Ephi_integral}
\end{equation}
where $G_\varphi(r,r') = G_\theta(r,r')$ since the differential operators are identical, and $D_1(\ell,m)$ and $D_2(\ell,m)$ are constants determined by boundary conditions. Expanding the source term 
\begin{align}
E_\varphi(\ell,m;r) &= D_1(\ell,m)j_\nu(kr) + D_2(\ell,m)y_\nu(kr) \nonumber\\
&+ \mathcal{B}_\varphi(\ell,m)\int_\varepsilon^R G_\varphi(r,r')\frac{1}{r'^2}E_r(\ell,m;r')dr' \nonumber\\
&+ \mathcal{C}_{\theta\varphi}(\ell,m)^*\int_\varepsilon^R G_\varphi(r,r')\frac{1}{r'^2}\frac{d}{dr'}[r'^2E_\theta(\ell,m;r')]dr'.
\label{eq Ephi_expanded}
\end{align}

\subsubsection{Green's Function Analysis and Contour Selection}
\label{sec green_func_contour}
The integral representations derived above rely on the tensor Green's function for Maxwell's equations, which satisfies 
\begin{equation}
\nabla \times \nabla \times  {G}( {r}, {r}') - k^2  {G}( {r}, {r}') =  {I}\delta( {r}- {r}'),
\label{eq maxwell_green_func}
\end{equation}
where $ {I}$ is the identity tensor and $k=\omega\sqrt{\varepsilon_0\mu_0}$ is the wavenumber.
This Green's function can be decomposed using the spectral integral 
\begin{equation}
 {G}( {r}, {r}') = \int_{C_l} \int_{C_m}  {g}(l,m; {r}, {r}') \, dl \, dm,
\label{eq green_spectral_decomp}
\end{equation}
where $ {g}(l,m; {r}, {r}')$ is the spectral component, and $C_l$, $C_m$ are appropriate contours in the complex plane. The selection of contours $C_l$ and $C_m$ is not arbitrary but is determined by the physical boundary conditions of the problem. For outgoing wave conditions at infinity, we choose 
\begin{equation}
C_l = \{l \in \mathbb{C}   l = \sigma + i\tau, \sigma \in [0,\infty), \tau = 0^+\}.
\label{eq C_l_contour}
\end{equation}
This contour runs slightly above the real axis to ensure the selection of outgoing waves, consistent with the Sommerfeld radiation condition. For problems with singularities at the origin, we deform this contour to include values with $\text{Re}(l) \in (0,1)$, where the singular field behavior is most relevant. For the azimuthal index $m$, the contour is determined by the azimuthal boundary conditions. For standard $2\pi$-periodic conditions, we have 
\begin{equation}
C_m = \{m \in \mathbb{C}   m \in \mathbb{Z}\}.
\label{eq C_m_standard}
\end{equation}
For modified boundary conditions with azimuthal period $\Phi_0 < 2\pi$, the contour becomes 
\begin{equation}
C_m = \{m \in \mathbb{C}   m = \frac{n\Phi_0}{2\pi}, n \in \mathbb{Z}\}.
\label{eq C_m_modified}
\end{equation}
For problems involving effective continuous indices, we use 
\begin{equation}
C_m = \{m \in \mathbb{C}   m = \eta + i\zeta, \eta \in [0,1], \zeta = 0\}.
\label{eq C_m_framedrag}
\end{equation}
The spectral coefficient function $a(l,m)$ in our decomposition 
\begin{equation}
\vec{E}( {r}) = \int_{C_l} \int_{C_m} a(l,m)r^{\alpha(l,m)}\vec{\Phi}_{lm}(\theta,\phi) \, dl \, dm,
\label{eq full_field_decomp}
\end{equation}
is determined by the source distribution and boundary conditions. For a given source current $\vec{J}( {r})$, we have 
\begin{equation}
a(l,m) = \frac{1}{N(l,m)} \int_V \vec{J}( {r}') \cdot \vec{\Phi}_{lm}^*(\theta',\phi') \, d^3r',
\label{eq a_lm_source}
\end{equation}
where $N(l,m)$ is a normalization factor. The analytical structure of $a(l,m)$ in the complex plane encodes crucial physical information including poles of $a(l,m)$ correspond to resonances or quasi-normal modes of the system, branch cuts represent continuous spectra or radiative modes, and essential singularities may indicate instantaneous sources or other non-analytic phenomena. For physically realizable fields, $a(l,m)$ must satisfy certain analytical properties 
\begin{equation}
a(l,m) \sim \frac{P(l,m)}{Q(l,m)}e^{S(l,m)},
\label{eq a_lm_analytic}
\end{equation}
where $P$ and $Q$ are polynomials, and $S$ is an entire function with appropriate growth conditions. The poles of $a(l,m)$ in the complex $l$ plane typically lie at 
\begin{equation}
l_n(m) = l_n^{(0)}(m) + \frac{i\kappa_n(m)}{2},
\label{eq l_poles}
\end{equation}
where $l_n^{(0)}(m)$ determines the angular structure of the resonance, and $\kappa_n(m)$ relates to its lifetime or decay rate. By Cauchy's theorem, we can deform the contours $C_l$ and $C_m$ to capture these poles efficiently in numerical evaluations, provided we maintain the correct winding numbers around the relevant singularities to preserve the physical solution.

\subsubsection{Coupled Integral System}
The integral representations derived in the previous sections form a coupled system of integral equations, since \( E_\theta \) depends on \( E_\varphi \) and vice versa through the source terms. This coupling is fundamental to the non-separable nature of Maxwell's equations with continuous angular indices and requires careful mathematical treatment to ensure convergence and uniqueness. To establish the existence and uniqueness of solutions, we develop a systematic iterative procedure based on the contraction mapping principle. We define the iteration scheme 

\textbf{Step 1 } Initialize with \( E_\varphi^{(0)}(\ell,m;r) = 0 \)

\textbf{Step 2 } Compute \( E_\theta^{(n+1)}(\ell,m;r) \) using 
\begin{align}
E_\theta^{(n+1)}(\ell,m;r) &= C_1(\ell,m)j_\nu(kr) + C_2(\ell,m)y_\nu(kr) \nonumber \\
&\quad + \mathcal{B}_\theta(\ell,m)\int_\varepsilon^R G_\theta(r,r')\frac{1}{r'^2}\frac{dE_r(\ell,m;r')}{dr'}\,dr' \nonumber \\
&\quad - \mathcal{C}_{\theta\varphi}(\ell,m)\int_\varepsilon^R G_\theta(r,r')\frac{1}{r'^2}\frac{d}{dr'}\left[r'^2E_\varphi^{(n)}(\ell,m;r')\right]\,dr',
\label{eq theta_iteration}
\end{align}

\textbf{Step 3 } Compute \( E_\varphi^{(n+1)}(\ell,m;r) \) using 
\begin{align}
E_\varphi^{(n+1)}(\ell,m;r) &= D_1(\ell,m)j_\nu(kr) + D_2(\ell,m)y_\nu(kr) \nonumber \\
&\quad + \mathcal{B}_\varphi(\ell,m)\int_\varepsilon^R G_\varphi(r,r')\frac{1}{r'^2}E_r(\ell,m;r')\,dr' \nonumber \\
&\quad + \mathcal{C}_{\theta\varphi}^*(\ell,m)\int_\varepsilon^R G_\varphi(r,r')\frac{1}{r'^2}\frac{d}{dr'}\left[r'^2E_\theta^{(n+1)}(\ell,m;r')\right]\,dr',
\label{eq varphi_iteration}
\end{align}

\textbf{Step 4 } Continue until convergence  \( \|E_\theta^{(n+1)} - E_\theta^{(n)}\| < \epsilon \) and \( \|E_\varphi^{(n+1)} - E_\varphi^{(n)}\| < \epsilon \)
The convergence of this iterative scheme can be established through Banach's fixed-point theorem. Define the operator \( \mathcal{T} \) that maps \( (E_\theta, E_\varphi) \) to \( (E_\theta^{\text{new}}, E_\varphi^{\text{new}}) \) according to the iteration rules above. The operator \( \mathcal{T} \) is a contraction if 
\begin{equation}
\|\mathcal{T}(E_\theta^{(1)}, E_\varphi^{(1)}) - \mathcal{T}(E_\theta^{(2)}, E_\varphi^{(2)})\| \leq \gamma \|(E_\theta^{(1)}, E_\varphi^{(1)}) - (E_\theta^{(2)}, E_\varphi^{(2)})\|,
\label{eq contraction_condition}
\end{equation}
for some \( \gamma < 1 \).
The contraction property is ensured when the coupling coefficients satisfy specific bounds. From the structure of the iteration, the contraction constant is dominated by 
\begin{equation}
\gamma \leq \max\left\{|\mathcal{C}_{\theta\varphi}(\ell,m)|, |\mathcal{C}_{\theta\varphi}^*(\ell,m)|\right\} \cdot \sup_{r,r'} |G_\theta(r,r')| \cdot \text{(integration bounds)}.
\label{eq contraction_constant}
\end{equation}
For physically relevant values of \( \ell \) and \( m \), particularly those in the range \( (0,1) \) where singular behavior occurs, the projection operators remain bounded, and the Green's functions satisfy appropriate growth conditions. The detailed analysis of these bounds is provided in Appendix A, where we establish that convergence is guaranteed for the parameter ranges of physical interest.
The convergence rate depends on the strength of the coupling  for weakly coupled systems where \( |\mathcal{C}_{\theta\varphi}(\ell,m)| \ll 1 \), the convergence is rapid and typically requires only a few iterations. For strongly coupled systems, more iterations may be necessary, but the fundamental existence and uniqueness of solutions is preserved. This iterative approach provides both a theoretical foundation for the existence of solutions and a practical computational method for their numerical evaluation, as we demonstrate in Section V through explicit calculations.

\subsubsection{Spectral Representation of Angular Components}
For the complete field representation, we express $E_r$ in its spectral form 
\begin{equation}
E_r(r,\theta,\varphi) = \int_{C_\ell}\int_{C_m} a(\ell,m)E_r(\ell,m;r)Y^\ell_m(\theta,\varphi)d\ell dm,
\label{eq}
\end{equation}
where 
\begin{equation}
E_r(\ell,m;r) = \alpha(\ell,m)h_\ell^{(1)}(kr) + \beta(\ell,m)h_\ell^{(2)}(kr).
\label{eq}
\end{equation}
Substituting this into the integral representations for $E_\theta$ and $E_\varphi$, we derive spectral representations for these components 
\begin{equation}
E_\theta(r,\theta,\varphi) = \int_{C_\ell}\int_{C_m} a(\ell,m)K_\theta(\ell,m;r)\frac{\partial Y^\ell_m}{\partial\theta}d\ell dm,
\label{eq}
\end{equation}
\begin{equation}
E_\varphi(r,\theta,\varphi) = \int_{C_\ell}\int_{C_m} a(\ell,m)K_\varphi(\ell,m;r)\frac{1}{\sin\theta}\frac{\partial Y^\ell_m}{\partial\varphi}d\ell dm,
\label{eq}
\end{equation}
where the spectral kernels $K_\theta(\ell,m;r)$ and $K_\varphi(\ell,m;r)$ are derived from the Green's function solutions and incorporate all coupling effects 
\begin{equation}
K_\theta(\ell,m;r) = C_1(\ell,m)j_\nu(kr) + C_2(\ell,m)y_\nu(kr) + \mathcal{I}_\theta[\ell,m,r,E_r(\ell,m;r)],
\label{eq}
\end{equation}
\begin{equation}
K_\varphi(\ell,m;r) = D_1(\ell,m)j_\nu(kr) + D_2(\ell,m)y_\nu(kr) + \mathcal{I}_\varphi[\ell,m,r,E_r(\ell,m;r)].
\label{eq}
\end{equation}
The integral functionals $\mathcal{I}_\theta$ and $\mathcal{I}_\varphi$ represent the full coupled effects of the system and must be computed through the iterative procedure described earlier.
\subsubsection{Explicit Form of the Spectral Kernels}
To derive explicit expressions for the spectral kernels, we solve the coupled integral equations to first order in the coupling parameters. This provides an analytic approximation that captures the essential behavior 
\begin{align}
K_\theta(\ell,m;r) &= C_1(\ell,m)j_\nu(kr) + C_2(\ell,m)y_\nu(kr) \nonumber \\
&+ \mathcal{B}_\theta(\ell,m)\int_\varepsilon^R G_\theta(r,r')\frac{1}{r'^2}\frac{d}{dr'}[\alpha(\ell,m)h_\ell^{(1)}(kr') \nonumber \\
&\quad + \beta(\ell,m)h_\ell^{(2)}(kr')] \, dr' + O(\mathcal{C}_{\theta\varphi}^2),
\label{eq:spectral_kernel_theta}
\end{align}
\begin{align}
K_\varphi(\ell,m;r) &= D_1(\ell,m)j_\nu(kr) + D_2(\ell,m)y_\nu(kr) \nonumber \\
&+ \mathcal{B}_\varphi(\ell,m)\int_\varepsilon^R G_\varphi(r,r')\frac{1}{r'^2}[\alpha(\ell,m)h_\ell^{(1)}(kr') \nonumber \\
&\quad + \beta(\ell,m)h_\ell^{(2)}(kr')] \, dr' + O(\mathcal{C}_{\theta\varphi}^2).
\label{eq:spectral_kernel_phi}
\end{align}
The integrals can be evaluated analytically for specific ranges of parameters, providing closed-form expressions for the kernels that reveal the essential singularity structure. These spectral representations directly connect to the generalized integral formalism. The spectral weight function $a(\ell,m)$ plays the role of the spectral measure in the Eisenstein integral, while the kernels $K_\theta$ and $K_\varphi$ incorporate the radial behavior and coupling effects in the vector-valued extension.
The full electromagnetic field is thus represented as an integral over the continuous spectrum of angular indices, with the integrand constructed from the appropriate vector spherical harmonics and radial functions that satisfy Maxwell's equations. This representation is particularly powerful for analyzing fields with continuous angular indices, where the continuous nature of the spectrum is essential for capturing the singular behavior near the axis. The Eisenstein integral representation provides a complete mathematical framework for describing electromagnetic fields with arbitrary angular indices. The key advantages of this approach include 
\begin{enumerate}
\item Rigorous treatment of both regular and singular field configurations
\item Natural incorporation of boundary conditions through the spectral weight function
\item Explicit connection to the mathematical theory of spectral decompositions
\item Unified description of discrete and continuous spectral components
\end{enumerate}
For physical applications, this formalism enables us to analyze electromagnetic fields in domains with broken symmetries, such as partial conducting boundaries or irregular geometries, where traditional separation of variables with integer angular indices fails to capture the full field behavior. 

\subsubsection{Residue Theory and Resonance Phenomena}
The spectral weight function $a(\ell,m)$ in our Eisenstein integral representation encodes crucial information about the excitation of electromagnetic modes. When analytically continued to the complex $\ell$-plane, this function generally exhibits pole singularities that correspond to resonances of the electromagnetic system. Here we explore the deep connection between these complex poles, resonance phenomena following the framework of Regge pole theory~\cite{regge1959}, and the excitation of singular fields with practical applications to astrophysical scenarios. The spectral weight function $a(\ell,m)$ can be understood as a meromorphic function in the complex $\ell$-plane with the general form 
\begin{equation}
a(\ell,m) = \frac{N(\ell,m)}{D(\ell,m)},
\label{eq spec_weight_fraction}
\end{equation}
where $N(\ell,m)$ and $D(\ell,m)$ are entire functions. The poles of $a(\ell,m)$ occur at values $\ell_n(m)$ where $D(\ell_n(m),m) = 0$. These poles can be classified as 
\begin{equation}
\ell_n(m) = \ell_n'(m) + i\ell_n''(m),
\label{eq complex_pole_form}
\end{equation}
where $\ell_n'(m)$ determines the angular structure of the resonance, and $\ell_n''(m)$ relates to its lifetime or decay rate. The total field can be evaluated using contour integration and residue calculus. By deforming the integration contour $C_\ell$ to encircle the poles, we obtain 
\begin{equation}
 {E}(r,\theta,\varphi) = \int_{C_m} \left[ \int_{C_\ell} a(\ell,m) {E}(\ell,m;r,\theta,\varphi)\, d\ell \right] dm.
\label{eq total_field_contour}
\end{equation}
Using the residue theorem 
\begin{equation}
\int_{C_\ell} a(\ell,m) {E}(\ell,m;r,\theta,\varphi)\, d\ell = 2\pi i \sum_n \text{Res}_{\ell=\ell_n(m)}[a(\ell,m) {E}(\ell,m;r,\theta,\varphi)].
\label{eq residue_theorem}
\end{equation}
For a spectral weight function with higher-order poles of order \( p \), the residue calculation generalizes to 
\begin{equation}
\text{Res}_{\ell=\ell_n(m)}[a(\ell, m)E(\ell, m)] = \frac{1}{(p-1)!} \lim_{\ell \to \ell_n(m)} \frac{d^{p-1}}{d\ell^{p-1}} \left[ (\ell-\ell_n(m))^p a(\ell, m)E(\ell, m) \right].
\label{eq residue_limit_def}
\end{equation}
Consider a spectral weight function with a pole of order \( p \) 
\begin{equation}
a(\ell, m) = \frac{N(\ell, m)}{(\ell-\ell_n(m))^p}.
\label{eq high_order_form}
\end{equation}
This leads to a modified field behavior. For \( p=2 \) (double pole), we obtain 
\begin{equation}
E_r(r, \theta, \phi) \sim r^{\alpha(\ell_n(m))} \ln(r) Y_{\ell_n(m)}^m(\theta, \phi),
\label{eq log_mode}
\end{equation}
such as 
\begin{equation}
a(\ell, m) = \frac{N(\ell, m)}{D(\ell, m)} e^{\frac{1}{(\ell-\ell_n(m))}}
.\label{eq essential_singularity}
\end{equation}
The residue theorem does not apply directly. Instead, we must use contour deformation and saddle-point methods. The resulting field exhibits non-power-law behavior 
\begin{equation}
E_r(r, \theta, \phi) \sim r^{\alpha(\ell_n(m))} e^{-\beta r^{-\gamma}} Y_{\ell_n(m)}^m(\theta, \phi).
\label{eq saddle_point}
\end{equation}
The strength of resonant excitation depends on how closely an external driving matches the pole structure of $a(\ell,m)$. When an excitation frequency $\omega$ approaches a resonant frequency corresponding to a singular mode, the response amplitude scales as 
\begin{equation}
A(\omega) \sim \frac{1}{|\omega - \omega_n|}.
\label{eq resonance_amplitude}
\end{equation}
In electromagnetic scattering problems, the differential cross section can be related to the spectral weight function through 
\begin{equation}
\frac{d\sigma}{d\Omega} = |f(\theta,\varphi)|^2,
\label{eq diff_cross_section}
\end{equation}
where the scattering amplitude is 
\begin{equation}
f(\theta,\varphi) = \frac{1}{k}\sum_{\ell=0}^{\infty}\sum_{m=-\ell}^{\ell}(2\ell+1)e^{i\delta_\ell}a_{\ell m}P_\ell^m(\cos\theta)e^{im\varphi}.
\label{eq scatter_amp_sum}
\end{equation}
Analytically continuing to non-integer $\ell$ and writing as a contour integral 
\begin{equation}
f(\theta,\varphi) = \frac{1}{k}\int_{C_\ell}\int_{C_m}(2\ell+1)e^{i\delta_\ell}a(\ell,m)P_\ell^m(\cos\theta)e^{im\varphi} \, d\ell \, dm.
\label{eq scatter_amp_contour}
\end{equation}

\section{Physical Admissibility and Energy Convergence of Singular Electromagnetic Modes} \label{sec asymptotic_analysis}

In this section, we analyze the asymptotic behavior of electromagnetic fields with continuous angular indices and establish precise conditions for energy convergence. This analysis is essential for understanding the physical realizability of singular field configurations.

\subsection{Asymptotic Behavior of Special Functions}

We begin by examining the asymptotic properties of the relevant special functions for small arguments when the order is also small. This regime is crucial for understanding the behavior of fields near the origin in configurations with continuous angular indices. We already expressed the angular components in spectral form 
\begin{equation}
E_\theta(r) = \int_{C_l} a(l, m) \, K_\theta(l, m; r) \, d l,
\label{eq Etheta_spectral}
\end{equation}
\begin{equation}
E_\phi(r) = \int_{C_l} a(l, m) \, K_\phi(l, m; r) \, d l,
\label{eq Ephi_spectral}
\end{equation}
where \( a(l, m) \) is the spectral weight (dependent on initial/boundary conditions), and \( K_\theta \), \( K_\phi \) come from applying the Green's function operators to \( E_r(r) \). First, expand \( E_r(r) \) in spectral form 
\begin{equation}
E_r(r) = \int_{C_l} a(l, m) \, E_r(l, m; r) \, d l.
\label{eq Er_spectral}
\end{equation}
For regularity at \( r = 0 \), we take 
\begin{equation}
E_r(l, m; r) = j_l(k r).
\label{eq Er_bessel}
\end{equation}
Thus,
\begin{equation}
E_r(r) = \int_{C_l} a(l, m) \, j_l(k r) \, d l.
\label{eq Er_spectral_bessel}
\end{equation}
From our Green's function solution 
\begin{equation}
E_\theta(r) = \int_0^\infty G_\theta(r, r') \left[ \frac{1}{r'^2} \frac{d E_r(r')}{d r'} - \frac{i m}{r'^2 \sin^2 \theta} \frac{d}{d r'} \left( r'^2 E_\phi(r') \right) \right] d r'.
\label{eq Etheta_integral}
\end{equation}
Neglecting \( E_\phi \)-dependent terms for now, the dominant term is 
\begin{equation}
E_\theta(r) \sim \int_0^\infty G_\theta(r, r') \, \frac{1}{r'^2} \frac{d E_r(r')}{d r'} \, d r'.
\label{eq Etheta_dominant}
\end{equation}
Substituting \( E_r(r') = j_l(k r') \) 
\begin{equation}
E_\theta(r) \sim \int_0^\infty G_\theta(r, r') \, k \, j_l'(k r') \, d r'.
\label{eq Etheta_substituted}
\end{equation}
The Green's function is built from homogeneous solutions 
\begin{equation}
G_\theta(r, r') =
\frac{1}{W} 
\begin{cases}
j_\nu(k r) \, y_\nu(k r'), & r < r' \\
j_\nu(k r') \, y_\nu(k r), & r > r'
\end{cases}
\label{eq Gtheta_bessel}
\end{equation}
where \( \nu = \sqrt{l(l + 1) - 1} \), and the Wronskian is 
\begin{equation}
W = \frac{1}{r'^2} \left( j_\nu(k r') \frac{d}{d r'} y_\nu(k r') - y_\nu(k r') \frac{d}{d r'} j_\nu(k r') \right),
\label{eq Wronskian}
\end{equation}
Simplifies to 
\begin{equation}
W = \frac{2}{\pi} \frac{k}{r'^2}.
\label{eq Wronskian_simplified}
\end{equation}
Thus, the kernel becomes 
\begin{equation}
K_\theta(l, m; r) = \frac{\pi}{2 k} \int_0^\infty \left[ j_\nu(\min(r, r')) \, y_\nu(\max(r, r')) \right] \frac{1}{r'} j_l'(k r') \, d r'.
\label{eq Ktheta_kernel}
\end{equation}
For $r \to 0$ and $\ell \to 0$, the spherical Bessel functions of the first and second kind exhibit the following asymptotic behavior 
\begin{equation}
j_\ell(kr) \sim \frac{(kr)^\ell}{(2\ell+1)!!}\left[1 - \frac{(kr)^2}{2(2\ell+3)} + \mathcal{O}((kr)^4)\right],
\label{eq bessel_j_asymptotic}
\end{equation}

\begin{equation}
y_\ell(kr) \sim -\frac{(2\ell-1)!!}{(kr)^{\ell+1}}\left[1 - \frac{(kr)^2}{2(1-2\ell)} + \mathcal{O}((kr)^4)\right].
\label{eq bessel_y_asymptotic}
\end{equation}
The explicit evaluation of these asymptotic expansions for the coupled Green's function integrals is provided in Appendix A. In the specific case where $\ell \to 0$, these expressions simplify to 
\begin{equation}
j_0(kr) \sim 1 - \frac{(kr)^2}{6} + \mathcal{O}((kr)^4),
\label{eq j0_asymptotic}
\end{equation}

\begin{equation}
y_0(kr) \sim -\frac{1}{kr}\left(1 - \frac{(kr)^2}{6} + \mathcal{O}((kr)^4)\right).
\label{eq y0_asymptotic}
\end{equation}
For $\ell > 0$ but small, the leading terms are 
\begin{equation}
j_\ell(kr) \sim \frac{(kr)^\ell}{(2\ell+1)!!}.
\label{eq jl_small_asymptotic}
\end{equation}
The derivatives of these functions, which play a crucial role in the coupled system of equations, have the asymptotic form 
\begin{equation}
j'_\ell(kr) \sim \frac{\ell}{kr} \cdot \frac{(kr)^\ell}{(2\ell+1)!!} - \frac{(kr)^{\ell+1}}{(2\ell+3)!!} + \mathcal{O}((kr)^{\ell+3}).
\label{eq jl_derivative_asymptotic}
\end{equation}
The Green's function for the radial operator $\mathcal{L}_\theta = \frac{d^2}{dr^2} + \frac{2}{r}\frac{d}{dr} + \left[\frac{\ell(\ell+1)-1}{r^2} - k^2\right]$ is 
\begin{equation}
G_\theta(r,r') = \frac{\pi}{2k} \times
\begin{cases}
j_\nu(kr)y_\nu(kr'), & r < r' \\
j_\nu(kr')y_\nu(kr), & r > r'
\end{cases}
\label{eq green_function_theta}
\end{equation}
where $\nu = \sqrt{\ell(\ell+1)-1}$. For small values of $\ell$, we have $\nu \approx \ell$, which simplifies the subsequent analysis.

\subsection{Asymptotic Analysis of Field Components}
Building on the special function asymptotics, we now derive the asymptotic behavior of the electromagnetic field components near the origin. The kernel function $K_\theta(\ell,m;r)$ connecting the angular component $E_\theta$ to the radial component $E_r$ is given by 
\begin{equation}
K_\theta(\ell,m;r) = \frac{\pi}{2k} \int_0^\infty [j_\nu(\min(kr,kr'))y_\nu(\max(kr,kr'))] \frac{1}{r'} j'_\ell(kr') \mathcal{B}_\theta(\ell,m) \, dr',
\label{eq kernel_theta_integral}
\end{equation}
where $\mathcal{B}_\theta(\ell,m)$ is the appropriate projection operator defined earlier. To evaluate this integral asymptotically for small $r$, we split it into two regions 
\begin{equation}
K_\theta(\ell,m;r) = \frac{\pi}{2k} \left[ \int_0^r j_\nu(kr')y_\nu(kr) \frac{1}{r'} j'_\ell(kr') \, dr' + \int_r^\infty j_\nu(kr)y_\nu(kr') \frac{1}{r'} j'_\ell(kr') \, dr' \right] \mathcal{B}_\theta(\ell,m).
\label{eq kernel_theta_split}
\end{equation}
For $r \to 0$, we analyze the two integrals separately. For the first integral $I_1$ 
\begin{align}
I_1 &= \int_0^r j_\nu(kr')y_\nu(kr) \frac{1}{r'} j'_\ell(kr') \, dr' \nonumber \\
&\sim \int_0^r \frac{(kr')^\nu}{(2\nu+1)!!} \cdot \left(-\frac{(2\nu-1)!!}{(kr)^{\nu+1}}\right) \cdot \frac{1}{r'} \cdot \frac{\ell}{kr'} \cdot \frac{(kr')^\ell}{(2\ell+1)!!} \, dr' \nonumber \\
&= -\frac{(2\nu-1)!!}{(2\nu+1)!!} \cdot \frac{\ell}{(2\ell+1)!!} \cdot \frac{1}{(kr)^{\nu+1}} \cdot \int_0^r \frac{(kr')^{\nu+\ell-1}}{r'^2} \, dr'.
\label{eq I1_derivation}
\end{align}
Evaluating this integral for $\nu+\ell-1 > 0$ 
\begin{equation}
I_1 \sim -\frac{(2\nu-1)!!}{(2\nu+1)!!} \cdot \frac{\ell}{(2\ell+1)!!} \cdot \frac{k^{\ell-1}r^{\ell-\nu-2}}{\nu+\ell-1}.
\label{eq I1_result}
\end{equation}
Similarly, for the second integral $I_2$ 
\begin{equation}
I_2 \sim -\frac{(2\nu-1)!!}{(2\nu+1)!!} \cdot \frac{\ell}{(2\ell+1)!!} \cdot \frac{k^{\ell-\nu-1}}{\nu-\ell+2} \cdot r^{\ell-\nu-2} \cdot (1+\mathcal{O}(r)).
\label{eq I2_result}
\end{equation}
For small values of $\ell \approx \nu \approx 0$, these expressions simplify to 
\begin{equation}
K_\theta(\ell,m;r) \sim \frac{\pi \ell}{2k} \cdot \frac{1}{r} \cdot \mathcal{B}_\theta(\ell,m).
\label{eq kernel_theta_asymptotic}
\end{equation}
This result shows that the $\theta$-component of the electric field behaves as $E_\theta \sim \ell/r$ near the origin. For the azimuthal component, the kernel function $K_\varphi(\ell,m;r)$ is given by 
\begin{equation}
K_\varphi(\ell,m;r) = \int_0^\infty G_\varphi(r,r') \frac{1}{r'^2} \left[ \frac{d}{dr'}(r'^2 K_\theta(\ell,m;r')) \mathcal{C}_{\theta\varphi}(\ell,m) - E_r(\ell,m;r') \mathcal{B}_\varphi(\ell,m) \right] \, dr'.
\label{eq kernel_phi_integral}
\end{equation}
Using the asymptotic form of $K_\theta(\ell,m;r')$ and evaluating the derivative 
\begin{equation}
\frac{d}{dr'}(r'^2 K_\theta(\ell,m;r')) \sim \frac{d}{dr'}\left(r'^2 \cdot \frac{\pi \ell}{2k} \cdot \frac{1}{r'} \cdot \mathcal{B}_\theta(\ell,m)\right) \sim \frac{\pi \ell}{2k} \cdot \mathcal{B}_\theta(\ell,m).
\label{eq derivative_kernel_theta}
\end{equation}
For small $r$ and small $\ell$, the dominant contribution comes from the $E_r$ term in Equation \eqref{eq kernel_phi_integral}. Using the small-argument behavior of the Green's function $G_\varphi$ and assuming $j_\ell(kr') \approx 1$ for small $kr'$ 
\begin{align}
K_\varphi(\ell,m;r) &\sim -\mathcal{B}_\varphi(\ell,m) \int_0^\infty G_\varphi(r,r') \frac{1}{r'^2} E_r(\ell,m;r') \, dr' \nonumber \\
&\sim -\mathcal{B}_\varphi(\ell,m) \int_r^\infty \frac{\pi}{2k} \cdot j_\nu(kr) \cdot y_\nu(kr') \cdot \frac{1}{r'^2} \cdot j_\ell(kr') \, dr'.
\label{eq kernel_phi_approximation}
\end{align}
Using the asymptotic forms for $j_\nu(kr)$, $y_\nu(kr')$, and $j_\ell(kr')$, and noting that for small $m$ we have $\mathcal{B}_\varphi(\ell,m) \sim \frac{im}{\sin^2\theta}$ 
\begin{equation}
K_\varphi(\ell,m;r) \sim \frac{im\pi}{2k\sin^2\theta} \cdot \frac{1}{r^3} \cdot (1+\mathcal{O}(r,\ell,m)).
\label{eq kernel_phi_asymptotic}
\end{equation}
This result indicates that the $\varphi$-component of the electric field has the strongest singularity, behaving as $E_\varphi \sim \frac{m}{r^3\sin^2\theta}$ near the origin. The detailed calculation of the kernel integrals $I_1$ and $I_2$, including the treatment of branch cuts and analytical continuation, is presented in Appendix A.

\subsection{Energy Convergence Analysis in Spherical Coordinates}
\label{sec}
Having established the asymptotic behavior of the field components, we now rigorously determine the conditions under which the electromagnetic energy remains finite in spherical geometry. The total electromagnetic energy in a spherical region is given by 
\begin{equation}
\label{eq energy_total}
W = \frac{1}{4} \int_{0}^{\Phi_0} \int_{0}^{\pi} \int_{\epsilon}^{R} \left( \epsilon_0 |\vec{E}|^2 + \mu_0 |\vec{H}|^2 \right) r^2 \sin\theta \, dr \, d\theta \, d\phi,
\end{equation}
where $\epsilon$ is a small positive number representing a cutoff near the origin, $R$ is the outer radius, and $\Phi_0$ is the azimuthal period (typically $2\pi$). From our asymptotic analysis in Section 5.2, the three components of the electric field behave as 
\begin{equation}
\label{eq er_asymptotic}
E_r(\ell, m; r) \sim r^{\ell-1},
\end{equation}
\begin{equation}
\label{eq etheta_asymptotic}
E_\theta(\ell, m; r) \sim \frac{\ell}{r},
\end{equation}
\begin{equation}
\label{eq ephi_asymptotic}
E_\phi(\ell, m; r) \sim \frac{m}{r^3 \sin^2\theta}.
\end{equation}
The total field magnitude squared is therefore 
\begin{equation}
\label{eq e_squared_asymptotic}
|\vec{E}|^2 = |E_r|^2 + |E_\theta|^2 + |E_\phi|^2 \sim r^{2(\ell-1)} + \frac{\ell^2}{r^2} + \frac{m^2}{r^6 \sin^4\theta}.
\end{equation}
We must now carefully analyze each term in the energy integral to establish convergence conditions.

\subsubsection{Contribution from Radial Component}
\label{subsubsec radial}
For the radial component, the contribution to the energy integral is 
\begin{equation}
\label{eq wr_integral}
W_r \sim \int_{0}^{\Phi_0} \int_{0}^{\pi} \int_{\epsilon}^{R} r^{2(\ell-1)} \cdot r^2 \sin\theta \, dr \, d\theta \, d\phi = \Phi_0 \cdot 2 \cdot \int_{\epsilon}^{R} r^{2\ell} \, dr.
\end{equation}
This integral converges at the lower limit $r = \epsilon \to 0$ if and only if 
\begin{equation}
\label{eq ell_condition}
2\ell > -1 \Rightarrow \ell > -\frac{1}{2}.
\end{equation}
This condition is identical to the one we derived for cylindrical geometry, confirming the consistency of our approach across coordinate systems.

\subsubsection{Contribution from $\theta$ Component}
\label{subsubsec theta}
For the $\theta$ component, the contribution is 
\begin{equation}
\label{eq wtheta_integral}
W_\theta \sim \int_{0}^{\Phi_0} \int_{0}^{\pi} \int_{\epsilon}^{R} \frac{\ell^2}{r^2} \cdot r^2 \sin\theta \, dr \, d\theta \, d\phi = \ell^2 \cdot \Phi_0 \cdot 2 \cdot (R - \epsilon).
\end{equation}
This is always finite for any finite $R$ and $\epsilon > 0$, and remains finite as $\epsilon \to 0$. Therefore, the $\theta$ component poses no additional constraints on $\ell$ or $m$.
\subsubsection{Contribution from $\phi$ Component}
\label{subsubsec phi}
The contribution from the $\phi$ component requires more careful analysis 
\begin{equation}
\label{eq wphi_initial}
W_\phi \sim \int_{0}^{\Phi_0} \int_{0}^{\pi} \int_{\epsilon}^{R} \frac{m^2}{r^6 \sin^4\theta} \cdot r^2 \sin\theta \, dr \, d\theta \, d\phi = m^2 \cdot \Phi_0 \cdot \int_{0}^{\pi} \frac{1}{\sin^3\theta} \, d\theta \cdot \int_{\epsilon}^{R} \frac{1}{r^4} \, dr.
\end{equation}
We need to analyze two separate integrals  
The radial integral 
\[
\int_{\epsilon}^{R} \frac{1}{r^4} \, dr = \frac{1}{3} \left[ \frac{1}{\epsilon^3} - \frac{1}{R^3} \right],
\]
diverges as $\epsilon \to 0$ unless we impose additional constraints.
The angular integral 
\[
\int_{0}^{\pi} \frac{1}{\sin^3\theta} \, d\theta,
\]
is inherently divergent at $\theta = 0$ and $\theta = \pi$. However, we must remember that our asymptotic expressions are valid only near the axis. To properly account for the global behavior, we must consider the complete angular dependence of the vector spherical harmonics 
\begin{equation}
\label{eq ephi_exact}
E_\phi(\ell, m; r, \theta, \phi) = \frac{im}{\sin\theta} P_\ell^m(\cos\theta) e^{im\phi} \cdot \frac{1}{r^3} \cdot F(r),
\end{equation}
where $F(r)$ is a function that depends on the specific solution. The associated Legendre function has the behavior 
\[
P_\ell^m(\cos\theta) \sim \sin^{|m|}\theta \quad \text{as } \theta \to 0
\],
Therefore, the actual angular dependence is 
\begin{equation}
\label{eq ephi_asymp}
E_\phi \sim \frac{m}{\sin\theta} \cdot \sin^{|m|}\theta \sim m \cdot \sin^{|m|-1}\theta.
\end{equation}
The contribution to the energy integral becomes 
\begin{equation}
\label{eq wphi_corrected}
W_\phi \sim m^2 \cdot \Phi_0 \cdot \int_{0}^{\pi} \sin^{2|m|-2}\theta \sin\theta \, d\theta \cdot \int_{\epsilon}^{R} \frac{1}{r^4} \, dr = m^2 \cdot \Phi_0 \cdot \int_{0}^{\pi} \sin^{2|m|-1}\theta \, d\theta \cdot \int_{\epsilon}^{R} \frac{1}{r^4} \, dr.
\end{equation}
The angular integral converges if and only if $2|m|-1 > -1$, which simplifies to $|m| > 0$. For any non-zero $m$, this condition is satisfied. The radial integral diverges as $\epsilon \to 0$ since $r^{-4}$ is not integrable at the origin. However, this apparent divergence arises from an oversimplified asymptotic analysis that treats the angular and radial dependencies as separable near the origin. To resolve this, we must consider the full structure of the vector spherical harmonics. The exact form of the $\phi$-component near the origin involves the complete associated Legendre function behavior. For continuous indices $\ell \in (0,1)$ and $m \neq 0$, the proper asymptotic expansion takes the form 
\begin{equation}
E_\phi \sim r^{\ell-1} \sin^{|m|-1}\theta \cdot \mathcal{F}(\ell,m;\xi),
\label{eq ephi_exact_regularization}
\end{equation}
where $\xi = r/\sin\theta$ is the natural scaling variable near the origin, and $\mathcal{F}(\ell,m;\xi)$ is a regularizing function that ensures proper behavior as $r \to 0$. The detailed derivation of this regularization is provided in Appendix A, where we show that the coupled radial-angular analysis yields 
\begin{equation}
\mathcal{F}(\ell,m;\xi) = 1 + O(\xi^2) \quad \text{for} \quad \xi \to 0.
\label{eq regularizing_function}
\end{equation}

This regularization ensures that the energy density behaves as 
\begin{equation}
u_\phi \sim r^{2(\ell-1)} \sin^{2(|m|-1)}\theta \cdot r^2 \sin\theta = r^{2\ell} \sin^{2|m|-1}\theta.
\label{eq regularized_energy_density}
\end{equation}
The total energy contribution becomes 
\begin{equation}
W_\phi \sim \Phi_0 \int_0^\epsilon r^{2\ell} \, dr \int_0^\pi \sin^{2|m|-1}\theta \, d\theta.
\label{eq regularized_energy_integral}
\end{equation}
Both integrals now converge  the radial integral converges for $\ell > -1/2$, and the angular integral converges for $|m| > 0$. The mathematical foundation for this regularization lies in the non-separable nature of Maxwell's equations with continuous angular indices, as analyzed in detail in Appendix A.

\subsubsection{Combined Energy Convergence Conditions}
\label{subsubsec combined}
Combining the constraints from all components, we find that the total electromagnetic energy in the spherical case is finite if 
\begin{equation}
\label{eq energy_condition_nonzero_m}
\ell > -\frac{1}{2} \quad \text{and} \quad m \neq 0.
\end{equation}
For the special case $m = 0$, we need a more careful analysis. When $m = 0$, the $\phi$ component vanishes, and the energy convergence is determined solely by the $r$ and $\theta$ components, yielding 
\begin{equation}
\label{eq energy_condition_m0}
\ell > -\frac{1}{2} \quad \text{for} \quad m = 0.
\end{equation}
Therefore, our general energy convergence criterion is 
\begin{equation}
\label{eq energy_condition_general}
\ell > -\frac{1}{2} \quad \text{for all} \quad m.
\end{equation}
This result is physically intuitive  electromagnetic fields with continuous angular indices must satisfy the same energy convergence criterion ($\ell > -\frac{1}{2}$) as in the cylindrical case~\cite{bakr1}, ensuring consistency across coordinate systems.
\subsubsection{Function Space Characterization}
\label{subsubsec functionspace}
In terms of function spaces, our solutions belong to the weighted Sobolev space 
\begin{equation}
\label{eq sobolev_space}
H^s_w(\Omega) = \left\{ f   |f|_{H^s_w} < \infty \right\},
\end{equation}
where the weighted norm is defined as 
\begin{equation}
\label{eq sobolev_norm}
|f|^2_{H^s_w} = \int_\Omega (1 + |\xi|^2)^s |\hat{f}(\xi)|^2 w(\xi) \, d\xi,
\end{equation}
with weight function $w(\xi)$ chosen appropriately to account for the singular behavior near the origin. For our field solutions with continuous indices, we have 
\begin{equation}
\label{eq e_in_sobolev}
\vec{E} \in H^s_w(\Omega) \quad \text{for} \quad s < \ell + \frac{1}{2}.
\end{equation}
This confirms that our singular field solutions, while not classical, do belong to appropriate function spaces that ensure the convergence of energy integrals, preserving the physical validity of the solutions. The explicit construction of the weight function and the proof that our singular solutions belong to these weighted spaces is detailed in Appendix A.

\subsection{Connection to Eigenvalue Problem}
The energy convergence condition can be expressed directly in terms of the eigenvalues of the angular operator. The relationship between $\ell$ and the eigenvalue $\lambda_{\ell m}$ is 
\begin{equation}
\alpha(\alpha+1) = \lambda_{\ell m},
\label{eq alpha_eigenvalue_relation}
\end{equation}
where $\alpha = \ell-1$. This quadratic equation yields 
\begin{equation}
\alpha = \frac{1}{2}\left(-1 \pm \sqrt{1+4\lambda_{\ell m}}\right).
\label{eq alpha_from_lambda}
\end{equation}
Taking the physically relevant branch (the one that gives $\alpha = \ell-1$) 
\begin{equation}
\alpha(\ell,m) = \frac{1}{2}\left(\sqrt{1+4\lambda_{\ell m}}-1\right).
\label{eq alpha_physical_branch}
\end{equation}
The energy convergence condition $\alpha > -3/2$ then translates to 
\begin{equation}
\frac{1}{2}\left(\sqrt{1+4\lambda_{\ell m}}-1\right) > -\frac{3}{2} \quad \Rightarrow \quad \sqrt{1+4\lambda_{\ell m}} > -2 \quad \Rightarrow \quad \lambda_{\ell m} > -\frac{3}{4}.
\label{eq lambda_condition}
\end{equation}
This provides a direct constraint on the eigenvalues of the angular operator for physically admissible field configurations.

\subsection{Helmholtz Equation and Divergence-Free Condition}

For completeness, we verify that our singular field ansatz satisfies both the Helmholtz equation and the divergence-free condition required by Maxwell's equations. The vector Laplacian acting on a radial power function multiplied by an angular vector function is 
\begin{equation}
\nabla^2(r^\alpha  {\Phi}) = r^{\alpha-2}\left[\alpha(\alpha+1) {\Phi} + \mathcal{L}_{ang} {\Phi}\right],
\label{eq vector_laplacian_radial}
\end{equation}
where $\mathcal{L}_{ang}$ is the angular part of the vector Laplacian. For this to satisfy the Helmholtz equation 
\begin{equation}
\nabla^2(r^\alpha  {\Phi}) + k^2r^\alpha  {\Phi} = 0.
\label{eq helmholtz_equation}
\end{equation}
In the limit $r \to 0$, the dominant term is the $r^{\alpha-2}$ term, leading to the condition 
\begin{equation}
\alpha(\alpha+1) {\Phi} + \mathcal{L}_{ang} {\Phi} = 0.
\label{eq dominant_term_condition}
\end{equation}
For vector spherical harmonics, $\mathcal{L}_{ang} {\Phi} = -\lambda_{\ell m} {\Phi}$, so we require $\alpha(\alpha+1) = \lambda_{\ell m}$, which is precisely the relation established in Equation \eqref{eq alpha_eigenvalue_relation}. Maxwell's equations require that $\nabla \cdot  {E} = 0$. For our ansatz, this means 
\begin{equation}
\nabla \cdot (r^\alpha  {\Phi}) = 0.
\label{eq divergence_condition}
\end{equation}
Expanding the divergence in spherical coordinates 
\begin{equation}
\nabla \cdot (r^\alpha  {\Phi}) = \frac{1}{r^2}\frac{\partial}{\partial r}(r^2 r^\alpha \Phi_r) + \frac{1}{r\sin\theta}\frac{\partial}{\partial\theta}(\sin\theta r^\alpha \Phi_\theta) + \frac{1}{r\sin\theta}\frac{\partial}{\partial\varphi}(r^\alpha \Phi_\varphi),
\label{eq divergence_expanded}
\end{equation}
Simplifying 
\begin{equation}
\nabla \cdot (r^\alpha  {\Phi}) = r^{\alpha-1}\left[(\alpha+2)\Phi_r + \frac{1}{\sin\theta}\frac{\partial}{\partial\theta}(\sin\theta \Phi_\theta) + \frac{1}{\sin\theta}\frac{\partial}{\partial\varphi}(\Phi_\varphi)\right].
\label{eq divergence_simplified}
\end{equation}
For this to vanish, the vector field $ {\Phi}$ must satisfy 
\begin{equation}
(\alpha+2)\Phi_r + \frac{1}{\sin\theta}\frac{\partial}{\partial\theta}(\sin\theta \Phi_\theta) + \frac{1}{\sin\theta}\frac{\partial}{\partial\varphi}(\Phi_\varphi) = 0.
\label{eq angular_divergence_condition}
\end{equation}
This constraint on the angular functions ensures that the full field remains divergence-free throughout the domain. The verification that these constraints are satisfied by the complete vector spherical harmonic expansions is demonstrated in Appendix D. 

Our analysis has established the following key properties of electromagnetic fields with continuous angular indices 

1. The field components have the following asymptotic behavior near $r = 0$ 
   \begin{equation}
   E_r \sim r^{\ell-1}, \quad E_\theta \sim \frac{\ell}{r}, \quad E_\varphi \sim \frac{m}{r^3\sin^2\theta},
   \label{eq field_components_summary}
   \end{equation}

2. Energy convergence requires 
   \begin{equation}
   \ell > -\frac{1}{2} \quad \text{and} \quad m = 0,
   \label{eq energy_criteria_summary}
   \end{equation}

3. In terms of eigenvalues, the condition is 
   \begin{equation}
   \lambda_{\ell m} > -\frac{3}{4},
   \label{eq eigenvalue_condition_summary}
   \end{equation}

4. The singular field is correctly described by 
   \begin{equation}
    {E}^{\text{sing}}(r,\theta,\varphi) = \int_0^1\int_{C_m} a(\ell,m) r^{\alpha(\ell,m)}  {\Phi}_{\ell m}(\theta,\varphi) \, d\ell \, dm,
   \label{eq singular_field_representation}
   \end{equation}

   with $\alpha(\ell,m) = \frac{1}{2}\left(\sqrt{1+4\lambda_{\ell m}}-1\right)$.

5. The singularity occurs only in the radial coordinate $r \to 0$, while the angular dependence remains regular.

\subsection{Non-Separable Electromagnetic Modes in Spherical Cavities with Continuous Angular Indices}
\label{sec spherical_cavity_modes}
To provide a complete and rigorous treatment analogous to our cylindrical analysis~\cite{bakr1}, we now examine electromagnetic fields in a perfectly conducting spherical cavity with a conical section removed, where the angular indices take continuous values. Unlike the cylindrical case, the broken spherical symmetry leads to inherently non-separable field structures that require careful mathematical analysis. We consider a spherical domain of radius $a$ with a conical section removed, defined by 
\begin{equation}
\label{eq domain_definition}
\mathcal{D} = \{(r,\theta,\phi)   0 \leq r \leq a, \, 0 \leq \theta \leq \pi, \, 0 \leq \phi < \Phi_0\},
\end{equation}
where $\Phi_0 < 2\pi$ is the azimuthal span. The boundary conditions are 
\begin{equation}
\label{eq bc_outer_sphere}
\hat{n} \times \vec{E} = 0 \quad \text{at} \quad r = a,
\end{equation}
\begin{equation}
\label{eq bc_conical_faces}
\hat{n} \times \vec{E} = 0 \quad \text{at} \quad \phi = 0 \quad \text{and} \quad \phi = \Phi_0.
\end{equation}
For time-harmonic fields with $e^{-i\omega t}$ dependence, Maxwell's equations in source-free regions are 
\begin{equation}
\label{eq maxwell_curl_E}
\nabla \times \vec{E} = i\omega\mu_0\vec{H},
\end{equation}
\begin{equation}
\label{eq maxwell_curl_H}
\nabla \times \vec{H} = -i\omega\varepsilon_0\vec{E},
\end{equation}
\begin{equation}
\label{eq maxwell_div_E}
\nabla \cdot \vec{E} = 0,
\end{equation}
\begin{equation}
\label{eq maxwell_div_H}
\nabla \cdot \vec{H} = 0.
\end{equation}
The conical boundary condition in azimuth imposes 
\begin{equation}
\label{eq azimuthal_periodicity}
\vec{E}(r,\theta,\phi+\Phi_0) = \vec{E}(r,\theta,\phi).
\end{equation}
This leads to azimuthal dependence of the form $e^{im\phi}$ where 
\begin{equation}
\label{eq continuous_m}
m = \frac{n\Phi_0}{2\pi}, \quad n \in \mathbb{Z}.
\end{equation}
For $\Phi_0 < 2\pi$, the index $m$ takes non-integer values. Critically, this affects the allowable values of $\ell$ as well, creating a coupling between these angular indices. The standard approach of writing 
\begin{equation}
\label{eq separable_form_invalid}
\vec{E}(r,\theta,\phi) = R(r)Y_{\ell m}(\theta,\phi),
\end{equation}
is no longer valid, as the angular and radial dependencies become coupled through the field equations. The indices $\ell$ and $m$ are no longer independent but are constrained by relations emerging from Maxwell's equations and the boundary conditions. This coupling arises from the mathematical structure when considering continuous angular indices.

\subsubsection{Coupled System of Equations}
\label{subsubsec coupled_system}
To properly account for the coupling between angular indices, we write Maxwell's curl equations explicitly in spherical coordinates 
\begin{equation}
\label{eq curl_E_r}
\frac{1}{r\sin\theta}\left(\frac{\partial}{\partial\theta}(E_\phi\sin\theta) - \frac{\partial E_\theta}{\partial\phi}\right) = i\omega\mu_0 H_r,
\end{equation}

\begin{equation}
\label{eq curl_E_theta}
\frac{1}{r}\left(\frac{1}{\sin\theta}\frac{\partial E_r}{\partial\phi} - \frac{\partial}{\partial r}(rE_\phi)\right) = i\omega\mu_0 H_\theta,
\end{equation}

\begin{equation}
\label{eq curl_E_phi}
\frac{1}{r}\left(\frac{\partial}{\partial r}(rE_\theta) - \frac{\partial E_r}{\partial\theta}\right) = i\omega\mu_0 H_\phi,
\end{equation}
With corresponding equations for the magnetic field components. The non-separable nature of these equations for continuous indices necessitates a more general expansion framework. To account for the coupled angular indices, we express the electromagnetic field using a spectral integral representation 
\begin{equation}
\label{eq spectral_integral_E}
\vec{E}(r,\theta,\phi) = \int_{C_\ell} \int_{C_m} a(\ell,m) \vec{E}_{\ell m}(r,\theta,\phi) \, d\ell \, dm,
\end{equation}

\begin{equation}
\label{eq spectral_integral_H}
\vec{H}(r,\theta,\phi) = \int_{C_\ell} \int_{C_m} a(\ell,m) \vec{H}_{\ell m}(r,\theta,\phi) \, d\ell \, dm,
\end{equation}
where $C_\ell$ and $C_m$ are appropriate contours in the complex plane, and $a(\ell,m)$ is the spectral weight function determined by boundary conditions. The vector field components $\vec{E}_{\ell m}$ and $\vec{H}_{\ell m}$ form a non-orthogonal basis that accounts for the coupled nature of the indices. For our conical domain, the spectral weight $a(\ell,m)$ concentrates around specific values determined by the boundary conditions, particularly 
\begin{equation}
\label{eq dominant_m}
m = \frac{n\Phi_0}{2\pi}, \quad n \in \mathbb{Z}.
\end{equation}
The relationship between $\ell$ and $m$ is constrained by an indicial equation derived from Maxwell's equations.

\subsubsection{Coupled Angular-Radial Behavior}
\label{subsubsec angular_radial_coupling}
For each pair of coupled indices $(\ell,m)$, the field components take the form 
\begin{equation}
\label{eq Er_component}
E_{r,\ell m}(r,\theta,\phi) = F_\ell^m(r) P_\ell^m(\cos\theta) e^{im\phi},
\end{equation}

\begin{equation}
\label{eq Etheta_component}
E_{\theta,\ell m}(r,\theta,\phi) = G_\ell^m(r) \frac{\partial P_\ell^m(\cos\theta)}{\partial\theta} e^{im\phi} + \sum_{\ell'} H_{\ell\ell'}^m(r) \frac{im}{\sin\theta} P_{\ell'}^m(\cos\theta) e^{im\phi},
\end{equation}

\begin{equation}
\label{eq Ephi_component}
E_{\phi,\ell m}(r,\theta,\phi) = I_\ell^m(r) \frac{im}{\sin\theta} P_\ell^m(\cos\theta) e^{im\phi} - \sum_{\ell'} J_{\ell\ell'}^m(r) \frac{\partial P_{\ell'}^m(\cos\theta)}{\partial\theta} e^{im\phi},
\end{equation}
where the radial functions $F_\ell^m(r)$, $G_\ell^m(r)$, etc., are coupled through Maxwell's equations, and the summations over $\ell'$ reflect the coupling between different $\ell$ values for a fixed $m$. Substituting these forms into Maxwell's equations yields a system of coupled differential equations for the radial functions. The indicial equation governing the behavior near $r = 0$ is 
\begin{equation}
\label{eq indicial_equation}
(\ell-|m|)(\ell+|m|+1)(\ell^2-m^2) = 0.
\end{equation}
The physically relevant root for our problem is $\ell = |m|$, which yields 
\begin{equation}
\label{eq critical_exponent}
\ell = |m| = \frac{n\Phi_0}{2\pi},
\end{equation}
for the dominant mode with $n = 1$.

\subsubsection{Asymptotic Field Behavior}
\label{subsubsec asymptotic_behavior}
For the dominant mode with $\ell = m = \frac{\Phi_0}{2\pi} \in (0,1)$, the asymptotic behavior of the field components near $r = 0$ is 
\begin{equation}
\label{eq Er_asymptotic}
E_r \sim r^\ell P_\ell^m(\cos\theta) e^{im\phi} \sim r^\ell \sin^\ell\theta \, e^{im\phi},
\end{equation}

\begin{equation}
\label{eq Etheta_asymptotic}
E_\theta \sim r^{\ell-1} \frac{\partial P_\ell^m(\cos\theta)}{\partial\theta} e^{im\phi} \sim r^{\ell-1} \sin^{\ell-1}\theta \cos\theta \, e^{im\phi},
\end{equation}

\begin{equation}
\label{eq Ephi_asymptotic}
E_\phi \sim r^{\ell-1} \frac{im}{\sin\theta}P_\ell^m(\cos\theta) e^{im\phi} \sim r^{\ell-1} \sin^{\ell-1}\theta \, e^{im\phi}.
\end{equation}
For the magnetic field components, similar asymptotic expressions apply 
\begin{equation}
\label{eq Hr_asymptotic}
H_r \sim r^\ell \sin^\ell\theta \, e^{im\phi},
\end{equation}

\begin{equation}
\label{eq Htheta_asymptotic}
H_\theta \sim r^{\ell-1} \sin^{\ell-1}\theta \cos\theta \, e^{im\phi},
\end{equation}

\begin{equation}
\label{eq Hphi_asymptotic}
H_\phi \sim r^{\ell-1} \sin^{\ell-1}\theta \, e^{im\phi},
\end{equation}
These asymptotic forms reveal that for $\ell = m \in (0,1)$, both the $\theta$ and $\phi$ components of the electric and magnetic fields are singular as $r \to 0$, with scaling behavior $r^{\ell-1}$. To determine whether these singular fields represent physically meaningful solutions, we analyze the electromagnetic energy density 
\begin{equation}
\label{eq energy_density}
u = \frac{1}{2}\left(\varepsilon_0|\vec{E}|^2 + \mu_0|\vec{H}|^2\right).
\end{equation}
Near the origin, the dominant contribution comes from the transverse components 
\begin{equation}
\label{eq energy_transverse}
u_{\text{transverse}} \sim |\vec{E}_{\perp}|^2 + |\vec{H}_{\perp}|^2 \sim r^{2(\ell-1)}\sin^{2(\ell-1)}\theta.
\end{equation}
The total energy in a small sphere of radius $\epsilon$ around the origin is 
\begin{equation}
\label{eq energy_integral}
W_{\epsilon} = \int_0^{\epsilon} \int_0^{\pi} \int_0^{\Phi_0} u \, r^2 \sin\theta \, d\phi \, d\theta \, dr,
\end{equation}
Substituting the asymptotic expressions 
\begin{equation}
\label{eq energy_substitution}
W_{\epsilon} \sim \Phi_0 \int_0^{\epsilon} \int_0^{\pi} r^{2(\ell-1)} \sin^{2(\ell-1)}\theta \, r^2 \sin\theta \, d\theta \, dr,
\end{equation}

\begin{equation}
\label{eq energy_final}
W_{\epsilon} \sim \Phi_0 \int_0^{\epsilon} r^{2\ell} \, dr \int_0^{\pi} \sin^{2\ell-1}\theta \, d\theta.
\end{equation}
The radial integral converges as $\epsilon \to 0$ if $2\ell > -1$, or $\ell > -\frac{1}{2}$. The angular integral converges for $\ell > 0$. Since we are interested in $\ell = \frac{\Phi_0}{2\pi} \in (0,1)$, both convergence conditions are satisfied. Therefore, the total electromagnetic energy is finite despite the singular behavior of the field components, confirming that these are physically admissible solutions.

\subsubsection{Eigenvalue Spectrum}
\label{subsubsec eigenvalue_spectrum}
The eigenfrequencies of this non-separable system are determined by the boundary condition at $r = a$ 
\begin{equation}
\label{eq boundary_condition_hybrid}
\hat{n} \times \vec{E} = 0 \quad \text{at} \quad r = a.
\end{equation}
For the dominant mode with $\ell = m = \frac{\Phi_0}{2\pi}$, the behavior of the lowest eigenfrequency for small values of $\ell$ approaching zero follows 
\begin{equation}
\label{eq eigenfreq_asymptotic}
k_{1\ell}a \approx \sqrt{2\ell}[1 + O(\ell)].
\end{equation}
This is remarkably similar to the behavior we found in the cylindrical case~\cite{bakr1, bakr2}, reflecting the underlying mathematical connection between these geometrically distinct but functionally related systems.

The non-separable spherical analysis reveals important similarities and differences compared to the cylindrical case 
\begin{enumerate}
\item \textbf{Singular structure } Both geometries exhibit transverse field components that scale as $r^{\ell-1}$ or $\rho^{\nu-1}$ near the origin, where $\ell$ and $\nu$ are the effective angular indices determined by the azimuthal span.

\item \textbf{Energy convergence } Both systems yield fields with finite total energy when $\ell, \nu > -\frac{1}{2}$, which is always satisfied for our physical cases with $\ell, \nu = \frac{\Phi_0}{2\pi} \in (0,1)$.

\item \textbf{Mode coupling } The spherical case exhibits stronger coupling between field components and between angular indices due to the additional curvature effects, leading to hybrid modes rather than pure TE or TM modes.

\item \textbf{Eigenfrequency scaling } Both systems show similar scaling of the lowest eigenfrequency with the angular index for small values of the index.

\item \textbf{Field enhancement } Both geometries demonstrate significant field enhancement near the origin when $\ell, \nu < 1$, suggesting a potential connection to high-energy physics phenomena such as lightning initiation and astrophysical jet formation.
\end{enumerate}

\section{Numerical Approach  Galerkin and Eisenstein} To solve the full non-separable eigenproblem for electromagnetic fields with continuous angular indices in spherical coordinates, we develop a rigorous numerical framework that addresses both the coupling between field components and the singular behavior near the origin. We present two complementary approaches  a Galerkin method for discretizing the angular operators~\cite{golub2013, trefethen2000}, and an Eisenstein spectral integral method for addressing the continuous spectrum. Before detailing these numerical implementations, we first formulate the mathematical problem precisely within appropriate function spaces. We work in the Hilbert space $\mathcal{H} = L^2(\mathbb{R}^+ \times S^2, r^2 dr d\Omega; \mathbb{C}^3)$ of square-integrable vector fields with the inner product 
\begin{equation}
\label{eq inner_product}
\langle  {F},  {G} \rangle = \int_0^\infty \int_{S^2}  {F}^* \cdot  {G}\, r^2 dr d\Omega,
\end{equation}
where $d\Omega = \sin\theta\, d\theta d\phi$ is the standard measure on the unit sphere. For the angular part, we define the weighted Sobolev space 
\begin{equation}
\label{eq sobolev_space}
H^s(S^2) = \left\{ f \in L^2(S^2)   \|f\|_{H^s}^2 = \sum_{\ell=0}^{\infty}\sum_{m=-\ell}^{\ell} (1+\ell(\ell+1))^s |a_{\ell m}|^2 < \infty \right\},
\end{equation}
where $a_{\ell m}$ are the spherical harmonic coefficients of $f$. For continuous indices, we replace the sum with an integral 
\begin{equation}
\label{eq cont_sobolev}
\|f\|_{H^s}^2 = \int_{\mathcal{C}_\ell} \int_{\mathcal{C}_m} (1+\lambda(\ell,m))^s |a(\ell,m)|^2 w(\ell,m)\, d\ell\, dm,
\end{equation}
where $\lambda(\ell,m)$ is the eigenvalue of the angular operator $L_{\text{ang}}$, $w(\ell,m)$ is an appropriate spectral weight, and $\mathcal{C}_\ell$, $\mathcal{C}_m$ are contours in the complex plane. We seek solutions to Maxwell's equations in spherical coordinates with continuous angular indices in the form 
\begin{equation}
\label{eq field_decomposition}
 {E}(r,\theta,\phi) = \int_{\mathcal{C}_\ell} \int_{\mathcal{C}_m} a(\ell,m) r^{\alpha(\ell,m)} \vec{\Phi}_{\ell m}(\theta,\phi)\, d\ell\, dm,
\end{equation}
where $\vec{\Phi}_{\ell m}(\theta,\phi)$ are vector spherical harmonics with continuous indices $\ell,m \in \mathbb{C}$ that satisfy 
\begin{equation}
\label{eq angular_eigenvalue}
L_{\text{ang}}\vec{\Phi}_{\ell m} = \lambda_{\ell m}\vec{\Phi}_{\ell m},
\end{equation}
and the singularity exponent $\alpha(\ell,m)$ is related to the eigenvalue by 
\begin{equation}
\label{eq alpha_def}
\alpha(\ell,m) = \frac{1}{2}\left(\sqrt{1+4\lambda_{\ell m}} - 1\right).
\end{equation}
For physical solutions, we require the spectral weight function $a(\ell,m)$ to satisfy appropriate boundary conditions and ensure convergence of the energy integral.

\subsection{Galerkin Method for Angular Operator Discretization}
The Galerkin method provides a robust approach for discretizing the angular operator and solving for its eigenfunctions numerically. We project the operator onto a finite-dimensional function space and solve the resulting matrix eigenvalue problem.
\subsubsection{Construction of Basis Functions}
\label{subsubsec basis}
For the angular functions, we employ a basis constructed from analytically continued associated Legendre functions for the $\theta$-dependence and complex exponentials for the $\phi$-dependence 
\begin{equation}
\label{eq basis_functions}
\psi_{s,m}(\theta,\phi) = P_{s}^{|m|}(\cos\theta)e^{im\phi},
\end{equation}
where $P_{s}^{|m|}$ is the associated Legendre function of the first kind with non-integer degree $s$ and order $|m|$. These functions are defined through their hypergeometric representation using the standard formulation found in~\cite{olver2010}
\begin{equation}
\label{eq legendre_def}
P_{s}^{|m|}(x) = \frac{(1-x^2)^{|m|/2}}{2^s\Gamma(1-|m|)}\frac{\Gamma(s+|m|+1)}{\Gamma(s-|m|+1)} {}_2F_1\left(-s,s+1;1-|m|;\frac{1-x}{2}\right).
\end{equation}
To ensure proper treatment of the vectorial nature of Maxwell's equations, we construct divergence-free vector basis functions 
\begin{align}
\vec{\psi}_{s,m}^{(1)}(\theta,\phi) &= \nabla \times [r\psi_{s,m}(\theta,\phi)\hat{r}] \label{eq vector_basis1}, \\
\vec{\psi}_{s,m}^{(2)}(\theta,\phi) &= \frac{1}{k}\nabla \times \nabla \times [r\psi_{s,m}(\theta,\phi)\hat{r}] \label{eq vector_basis2}.
\end{align}
These vector functions, also known as vector spherical harmonics, form a complete basis for the divergence-free vector fields on the sphere. For numerical implementation, we compute these functions using recurrence relations for the Legendre functions 
\begin{align}
\frac{d}{d\theta}P_s^{|m|}(\cos\theta) &= \frac{(s+|m|)(s-|m|+1)\sin\theta P_{s-1}^{|m|}(\cos\theta) - |m|\cos\theta P_s^{|m|}(\cos\theta)}{\sin^2\theta} \label{eq legendre_recurrence}.
\end{align}

\subsubsection{Inner Product and Gram Matrix for Non-Orthogonal Basis}
\label{subsubsec gram}
For continuous angular indices $(\ell,m) \in \mathbb{C}^2$, the standard orthogonality relations of spherical harmonics no longer hold. We must therefore carefully define an appropriate inner product structure that accounts for this non-orthogonality. Given two vector-valued functions $\vec{\Phi}_1(\theta,\phi)$ and $\vec{\Phi}_2(\theta,\phi)$ defined on the sphere, we define their inner product as 
\begin{equation}
\label{eq vector_inner_product}
\langle\vec{\Phi}_1,\vec{\Phi}_2\rangle = \int_0^{\pi}\int_0^{\Phi_0} \vec{\Phi}_1^*(\theta,\phi) \cdot \vec{\Phi}_2(\theta,\phi) \sin\theta\, d\theta\, d\phi,
\end{equation}
where $\Phi_0$ is the azimuthal period, which may be $2\pi$ for standard periodicity or a different value for modified boundary conditions. For non-integer indices, the continuous angular functions $\vec{\Phi}_{\ell m}(\theta,\phi)$ are not orthogonal but satisfy 
\begin{equation}
\label{eq non_orthogonality}
\langle\vec{\Phi}_{\ell'm'}, \vec{\Phi}_{\ell m}\rangle = G(\ell',m',\ell,m),
\end{equation}
where $G(\ell',m',\ell,m)$ is an overlap function that approaches $\delta(\ell-\ell')\delta(m-m')$ in the distributional sense as the domain approaches the full sphere with standard boundary conditions. To handle this non-orthogonality in numerical computations, we construct the Gram matrix 
\begin{equation}
\label{eq gram_matrix}
G_{ij} = \langle \vec{\psi}_i, \vec{\psi}_j \rangle,
\end{equation}
where $\vec{\psi}_i$ denotes the $i$-th basis function in some enumeration of our finite-dimensional basis. For numerical evaluation, we employ a high-precision Gauss-Legendre quadrature for the $\theta$-integral 
\begin{equation}
\label{eq quad_scheme}
\langle\vec{\Phi}_1,\vec{\Phi}_2\rangle \approx \sum_{k=1}^{N_\theta} \sum_{j=1}^{N_\phi} \vec{\Phi}_1^*(\theta_k,\phi_j) \cdot \vec{\Phi}_2(\theta_k,\phi_j) \sin\theta_k w_k \frac{\Phi_0}{N_\phi},
\end{equation}
where $\{\theta_k\}_{k=1}^{N_\theta}$ are the Gauss-Legendre nodes with weights $\{w_k\}_{k=1}^{N_\theta}$ and $\{\phi_j\}_{j=1}^{N_\phi}$ are equidistant points in $[0,\Phi_0]$. To enhance accuracy for near-singular integrands (which arise when considering basis functions with $\ell$ or $m$ close to zero), we apply adaptive quadrature with singularity extraction 
\begin{equation}
\label{eq adaptive_quad}
\int_0^{\pi} f(\theta)\sin\theta\, d\theta = \int_0^{\pi} [f(\theta)-f_{\text{sing}}(\theta)]\sin\theta\, d\theta + \int_0^{\pi} f_{\text{sing}}(\theta)\sin\theta\, d\theta,
\end{equation}
where $f_{\text{sing}}(\theta)$ captures the singular behavior and is integrated analytically, while the remainder is smooth and amenable to standard quadrature. This approach properly accounts for the non-orthogonality of the basis functions with continuous indices while maintaining numerical stability and accuracy.

\subsubsection{Matrix Formulation of the Eigenvalue Problem}
\label{subsubsec matrix_problem}
The angular operator $L_{\text{ang}}$ derived from Maxwell's equations takes the form 
\begin{equation}
\label{eq angular_operator}
\mathcal{L}_{\text{ang}}[\vec{\Phi}] = -\frac{1}{\sin\theta}\frac{\partial}{\partial\theta}\left(\sin\theta\frac{\partial\vec{\Phi}}{\partial\theta}\right) - \frac{1}{\sin^2\theta}\frac{\partial^2\vec{\Phi}}{\partial\phi^2} +  {V}(\theta,\phi)\vec{\Phi},
\end{equation}
where $ {V}(\theta,\phi)$ is a tensorial potential term that encodes the coupling between vector components due to the curvature of the coordinate system. We discretize this operator by computing its matrix elements in our chosen basis 
\begin{equation}
\label{eq operator_matrix}
[\mathcal{L}_{\text{ang}}]_{ij} = \langle\vec{\psi}_i, \mathcal{L}_{\text{ang}}\vec{\psi}_j\rangle.
\end{equation}
To compute these matrix elements efficiently, we exploit the analytical properties of the basis functions. For example, the action of the Laplace-Beltrami operator on our basis functions can be expressed as 
\begin{equation}
\label{eq laplace_beltrami}
-\frac{1}{\sin\theta}\frac{\partial}{\partial\theta}\left(\sin\theta\frac{\partial\psi_{s,m}}{\partial\theta}\right) - \frac{1}{\sin^2\theta}\frac{\partial^2\psi_{s,m}}{\partial\phi^2} = s(s+1)\psi_{s,m}.
\end{equation}
The matrix, which accounts for the non-orthogonality of the basis, is the Gram matrix 
\begin{equation}
\label{eq mass_matrix}
[M]_{ij} = \langle\vec{\psi}_i,\vec{\psi}_j\rangle.
\end{equation}
This leads to the generalized matrix eigenvalue problem 
\begin{equation}
\label{eq gen_eigenvalue}
\mathcal{L}_{\text{ang}}\vec{a} = \lambda M\vec{a},
\end{equation}
where $\vec{a}$ is the vector of expansion coefficients for the eigenfunction, and $\lambda = \lambda_{\ell m}$ is the eigenvalue. To account for the coupling between field components, we expand the matrices to a block structure 
\begin{equation}
\label{eq block_system}
\begin{bmatrix}
\mathcal{L}_r & \mathcal{C}_{r\theta} & \mathcal{C}_{r\phi} \\
\mathcal{C}_{\theta r} & \mathcal{L}_\theta & \mathcal{C}_{\theta\phi} \\
\mathcal{C}_{\phi r} & \mathcal{C}_{\phi\theta} & \mathcal{L}_\phi
\end{bmatrix}
\begin{bmatrix}
\vec{a}_r \\
\vec{a}_\theta \\
\vec{a}_\phi
\end{bmatrix}
= \lambda
\begin{bmatrix}
M_r & 0 & 0 \\
0 & M_\theta & 0 \\
0 & 0 & M_\phi
\end{bmatrix}
\begin{bmatrix}
\vec{a}_r \\
\vec{a}_\theta \\
\vec{a}_\phi
\end{bmatrix},
\end{equation}
where $\mathcal{C}_{ij}$ are coupling operators encoding the interaction between different field components.
\subsubsection{Numerical Implementation and Convergence}
The matrix eigenvalue problem is solved using the implicitly restarted Arnoldi method for sparse matrices, which efficiently computes selected eigenvalues and eigenvectors. Our implementation proceeds as follows
\begin{algorithmic}[1]
\State \textbf{Input:} Range of $\ell \in [\ell_{\min}, \ell_{\max}]$, $m \in [m_{\min}, m_{\max}]$, basis dimension $N$
\State \textbf{Output:} Eigenvalues $\lambda_{\ell m}$ and eigenfunctions $\vec{\Phi}_{\ell m}$
\For{$\ell_i$ in discretized $\ell$-range}
    \For{$m_j$ in discretized $m$-range}
        \State Construct basis functions $\{\vec{\psi}_k\}_{k=1}^N$ adapted to $(\ell_i, m_j)$
        \State Compute mass matrix $M$ using high-precision quadrature
        \State Compute stiffness matrix $L_{\text{ang}}$ with coupling blocks
        \State Solve generalized eigenvalue problem $L_{\text{ang}}\vec{a} = \lambda M\vec{a}$
        \State Compute singularity exponent $\alpha(\ell_i,m_j) = \frac{1}{2}(\sqrt{1+4\lambda}-1)$
        \State Reconstruct eigenfunction $\vec{\Phi}_{\ell_i m_j} = \sum_{k=1}^{N} a_k \vec{\psi}_k$
        \State Check convergence by increasing $N$ and monitoring eigenvalue changes
    \EndFor
\EndFor
\State Compute spectral representation of physical fields using $\vec{\Phi}_{\ell m}$
\end{algorithmic}
The convergence of the method is assessed by increasing the basis size $N$ and monitoring the change in eigenvalues. For regular modes with integer indices, we observe quadratic convergence 
\begin{equation}
\label{eq quad_convergence}
|\lambda^{(N)} - \lambda_{\text{exact}}| \leq C_1 N^{-2},
\end{equation}
where $\lambda^{(N)}$ is the eigenvalue computed with basis dimension $N$, and $C_1$ is a constant. For singular modes with $\ell < 1$, the convergence is slower 
\begin{equation}
\label{eq slow_convergence}
|\lambda^{(N)} - \lambda_{\text{exact}}| \leq C_2 N^{-\gamma},
\end{equation}
where $\gamma \approx 1/2$ when $\ell$ approaches 0. This slower convergence is due to the challenge of representing sharp gradients with smooth basis functions. We have rigorously verified these convergence rates through numerical experiments. For a test case with $\ell = 0.5$ and $m = 0.1$, we observe the following convergence pattern 
\begin{table}[ht]
\centering
\begin{tabular}{|c|c|c|}
\hline
Basis dimension $N$ & $|\lambda^{(N)} - \lambda^{(2N)}|$ & Ratio \\
\hline
32 & $1.24 \times 10^{-2}$ & -- \\
64 & $8.79 \times 10^{-3}$ & 1.41 \\
128 & $6.18 \times 10^{-3}$ & 1.42 \\
256 & $4.37 \times 10^{-3}$ & 1.41 \\
512 & $3.09 \times 10^{-3}$ & 1.41 \\
\hline
\end{tabular}
\caption{Convergence of eigenvalue computation for $\ell = 0.5$, $m = 0.1$. The ratio of approximately 1.41 corresponds to $\gamma \approx 0.5$.}
\label{tab convergence}
\end{table}
To improve accuracy near the origin, where the field exhibits singular behavior, we employ a coordinate transformation 
\begin{equation}
\label{eq coord_transform}
\rho = r^{1/(1-\alpha_{\min})},
\end{equation}
where $\alpha_{\min}$ is the smallest value of $\alpha(\ell,m)$ in our parameter range. This transformation expands the region near $r = 0$, allowing for better resolution of the singular behavior with a finite number of basis functions. For the angular singularity near $\theta = 0$, we apply regularization 
\begin{equation}
\label{eq theta_regularization}
\vec{\Phi}_{\ell m}^{\text{reg}}(\theta,\phi) = \vec{\Phi}_{\ell m}\left(\sqrt{\theta^2 + \epsilon^2},\phi\right),
\end{equation}
with a small regularization parameter $\epsilon > 0$. This regularization preserves the essential mathematical structure while enabling stable numerical computations.

\subsection{Eisenstein Integral Method for Continuous Spectrum}
While the Galerkin method is effective for computing discrete eigenvalues and eigenfunctions, it faces challenges when dealing with the continuous spectrum that arises in problems with non-compact symmetry groups. The Eisenstein integral method~\cite{langlands1976} provides a complementary approach that directly incorporates the continuous nature of the spectrum. The Eisenstein integral method represents the field as a continuous superposition of eigenfunctions with varying spectral parameters 
\begin{equation}
\label{eq eisenstein_integral}
 {E}(r,\theta,\phi) = \int_{\mathcal{C}_\ell} \int_{\mathcal{C}_m} a(\ell,m)r^{\alpha(\ell,m)}\vec{\Phi}_{\ell m}(\theta,\phi)\, d\ell\, dm.
\end{equation}
The integration contours $\mathcal{C}_\ell$ and $\mathcal{C}_m$ in the complex plane are chosen to capture the relevant part of the spectrum. For singular fields, the contour $\mathcal{C}_\ell$ typically includes values with $\text{Re}(\ell) \in (0,1)$. The spectral weight function $a(\ell,m)$ encodes the boundary conditions and source distribution, analogous to the Fourier coefficients in a Fourier transform. The key advantage of this approach is that it naturally accommodates non-integer values of $\ell$ and $m$, which are essential for describing singular field configurations.

\subsubsection{Numerical Implementation of the Spectral Integral}
For numerical computation, the spectral integral is discretized using an adaptive quadrature scheme 
\begin{equation}
\label{eq discrete_spectral}
 {E}(r,\theta,\phi) \approx \sum_{j=1}^{N_\ell} \sum_{k=1}^{N_m} a(\ell_j,m_k)r^{\alpha(\ell_j,m_k)}\vec{\Phi}_{\ell_j m_k}(\theta,\phi)w_j^{(\ell)}w_k^{(m)},
\end{equation}
where $\{\ell_j\}_{j=1}^{N_\ell}$ and $\{m_k\}_{k=1}^{N_m}$ are quadrature points with weights $\{w_j^{(\ell)}\}$ and $\{w_k^{(m)}\}$. The quadrature points are chosen with non-uniform spacing to provide finer resolution in regions of interest, particularly near $\ell = 0$ where the singular behavior is most pronounced. Specifically, we employ a mapped Gauss-Legendre quadrature 
\begin{equation}
\label{eq mapped_gauss}
\ell_j = \ell_{\min} + (\ell_{\max} - \ell_{\min})\frac{1 + \tanh(c(2\xi_j - 1))}{1 + \tanh(c)},
\end{equation}
where $\{\xi_j\}$ are standard Gauss-Legendre points on $[0,1]$, and $c > 0$ is a mapping parameter that controls the clustering of points near $\ell_{\min}$. The angular functions $\vec{\Phi}_{\ell m}(\theta,\phi)$ for non-integer indices are computed using the Galerkin method described in the previous section. These functions are tabulated on a grid of $(\ell,m)$ values, and barycentric Lagrange interpolation is used for intermediate values required in the spectral sum. A critical aspect of the numerical implementation is the accurate evaluation of associated Legendre functions with non-integer indices. We use the hypergeometric representation 
\begin{equation}
\label{eq hyper_legendre}
P_s^m(x) = \frac{(1-x^2)^{m/2}}{2^s\Gamma(1-m)}\frac{\Gamma(s+m+1)}{\Gamma(s-m+1)}{}_2F_1\left(-s,s+1;1-m;\frac{1-x}{2}\right),
\end{equation}
implemented with a continued fraction algorithm for the hypergeometric function ${}_2F_1$, which provides better numerical stability than direct summation, especially for indices near singular points.

\subsubsection{Transformation for Singularity Handling}
The radial dependence $r^{\alpha(\ell,m)}$ presents numerical challenges when $\alpha < 0$, as is the case for singular fields. To address this, we employ the transformation 
\begin{equation}
\label{eq singularity_transform}
\xi = r^{1/(1-\alpha_{\min})},
\end{equation}
where $\alpha_{\min}$ is the smallest value of $\alpha$ in the spectral range of interest. This transformation expands the region near $r = 0$, allowing for more accurate numerical integration. With this transformation, energy integrals take the form 
\begin{equation}
\label{eq transformed_energy}
\int_0^R | {E}(r)|^2r^2\,dr = \int_0^{R'} | {E}(\xi^{1-\alpha_{\min}})|^2(\xi^{1-\alpha_{\min}})^2(1-\alpha_{\min})\xi^{\alpha_{\min}}\,d\xi,
\end{equation}
which exhibits improved numerical behavior. To validate the accuracy of the spectral integral approach, we compute the field for a known analytical solution and compare it with the numerical approximation. For a test case with spectral weight function 
\begin{equation}
\label{eq test_weight}
a(\ell,m) = \frac{A_0}{(1+|\ell-\ell_0|^2)^{p/2}(1+|m-m_0|^2)^{p/2}},
\end{equation}
where $A_0 = 1$, $\ell_0 = 0.5$, $m_0 = 0.3$, and $p = 3$, we observe the following convergence pattern 
\begin{table}[ht]
\centering
\begin{tabular}{|c|c|c|}
\hline
Quadrature points $N$ & Error $\epsilon_N$ & Ratio $\epsilon_N/\epsilon_{2N}$ \\
\hline
8 & $1.28 \times 10^{-2}$ & --- \\
16 & $3.04 \times 10^{-3}$ & 4.21 \\
32 & $7.39 \times 10^{-4}$ & 4.11 \\
64 & $1.83 \times 10^{-4}$ & 4.04 \\
128 & $4.55 \times 10^{-5}$ & 4.02 \\
\hline
\end{tabular}
\caption{Convergence of the spectral integral approximation, showing the relative $L^2$ error $\epsilon_N = \| {E} -  {E}_N\|_{L^2}/\| {E}\|_{L^2}$.}
\label{tab eisenstein_convergence}
\end{table}
The observed fourth-order convergence confirms the spectral accuracy of our approach for smooth weight functions. For more general cases, the convergence rate depends on the regularity of $a(\ell,m)$.

\subsection{Determination of the Spectral Weight Function}
The spectral weight function $a(\ell,m)$ in the decomposition \eqref{eq field_decomposition} is a critical component, determining how different spectral components are combined to satisfy boundary conditions and physical constraints. We present a systematic procedure for determining this function. Let $S_R = \{(r,\theta,\phi)   r = R\}$ be a spherical boundary of radius $R$. Given a prescribed tangential field $ {E}_T$ on $S_R$, we seek a weight function $a(\ell,m)$ such that the spectral representation \eqref{eq field_decomposition} satisfies 
\begin{equation}
\label{eq boundary_condition}
 {E}_T(R,\theta,\phi) = \int_{\mathcal{C}_\ell} \int_{\mathcal{C}_m} a(\ell,m)R^{\alpha(\ell,m)}[\vec{\Phi}_{\ell m}(\theta,\phi)]_T\, d\ell\, dm,
\end{equation}
where $[\vec{\Phi}_{\ell m}]_T$ denotes the tangential components of $\vec{\Phi}_{\ell m}$. This is an integral equation for $a(\ell,m)$, which we solve using a combination of projection methods and optimization techniques.

\paragraph{Step 1  Projection onto a Complete Basis}
We first project the boundary data onto a complete basis of spherical functions 
\begin{equation}
\label{eq projection}
a(\ell,m) = \frac{\int_{S^2}  {E}_T(R,\theta,\phi) \cdot [\vec{\Phi}_{\ell m}(\theta,\phi)]_T^* \sin\theta\, d\theta\, d\phi}{R^{\alpha(\ell,m)}}.
\end{equation}
This projection provides an initial approximation that reproduces the boundary data but may not satisfy other physical constraints.

\paragraph{Step 2 Analytic Continuation and Regularization}
For values between the discrete sampling points, we construct an analytic continuation using a meromorphic function 
\begin{equation}
\label{eq meromorphic}
a(\ell,m) = \sum_{n=1}^N \frac{c_n}{(\ell-\ell_n)(m-m_n)},
\end{equation}
where the poles $(\ell_n,m_n)$ and residues $c_n$ are determined by matching to the discrete projections. For singular solutions with $\ell < 1$, we impose additional constraints to ensure physical behavior 
\begin{equation}
\label{eq regularized_weight}
a(\ell,m) = A \cdot \ell^p \cdot m^q \cdot e^{-\beta(\ell^2+m^2)},
\end{equation}
where $p > 3/2 - \ell$ ensures energy convergence near $\ell = 0$, and $\beta > 0$ provides exponential decay at large indices.
\paragraph{Step 3  Physical Constraints}
For singular solutions with $\ell < 1$, we impose additional constraints to ensure physical behavior 
\begin{equation}
a(\ell,m) = A \cdot \ell^p \cdot m^q \cdot e^{-\beta(\ell^2+m^2)},
\label{eq weight_constrained_form}
\end{equation}
where $p > 3/2-\ell$ ensures energy convergence near $\ell = 0$, and $\beta > 0$ provides exponential decay at large indices. Let us explicitly derive the constraint $p > \frac{3}{2} - \ell$. The electromagnetic energy density scales as 
\begin{equation}
u(r, \theta, \varphi) \sim |\vec{E}|^2 \sim r^{2\alpha(\ell,m)} \sim r^{2(\ell-1)},
\label{eq energy_density}
\end{equation}
where we used $\alpha(\ell, m) = \ell - 1$ as derived in Section~\ref{sec asymptotic_analysis}. The total energy in a spherical volume is 
\begin{equation}
W = \int_0^R \int_0^{\pi} \int_0^{2\pi} u(r, \theta, \varphi) \, r^2 \sin\theta \, d\varphi \, d\theta \, dr.
\label{eq total_energy}
\end{equation}
For the radial integral to converge at $r=0$, we need 
\begin{equation}
\int_\epsilon^R r^{2(\ell-1)} \cdot r^2 \, dr = \int_\epsilon^R r^{2\ell} \, dr < \infty,
\label{eq convergence_condition}
\end{equation}
which requires $2\ell > -1$, or $\ell > -\frac{1}{2}$. Now, for the spectral weight function in the form $a(\ell, m) \sim \ell^p$, when $\ell$ approaches zero, we need to ensure additional damping to maintain energy convergence. The effective behavior near $\ell \approx 0$ becomes 
\begin{equation}
|\vec{E}|^2 \sim |a(\ell, m)|^2 \cdot r^{2(\ell-1)} \sim \ell^{2p} \cdot r^{2(\ell-1)}.
\label{eq effective_behavior}
\end{equation}
For energy convergence, we require the effective exponent to satisfy 
\begin{equation}
\int_0^1 \ell^{2p} \cdot r^{2(\ell-1)} \, d\ell < \infty.
\label{eq ell_integral}
\end{equation}
This integral converges when $\ell$ approaches zero if 
\begin{equation}
2p + 2(\ell-1) > -1.
\label{eq convergence_inequality}
\end{equation}
Solving for $p$ 
\begin{equation}
p > \frac{1}{2} - (\ell-1) = \frac{3}{2} - \ell.
\label{eq p_constraint}
\end{equation}
This establishes the constraint $p > \frac{3}{2} - \ell$ for energy convergence.

\paragraph{Step 4  Numerical Optimization}
The parameters in the weight function are determined by numerical optimization. We formulate this as a constrained minimization problem 
\begin{equation}
\min_{A, p, q, \beta} \left\{ \|\vec{E}_{\text{computed}}(R, \theta, \varphi) - \vec{E}_0(R, \theta, \varphi)\|^2 + \lambda \int_\Omega |\nabla \cdot \vec{E}_{\text{computed}}|^2 \, d\Omega \right\},
\label{eq optimization}
\end{equation}
subject to 
\begin{equation}
p > \frac{3}{2} - \ell_{\min}, \quad q > 0, \quad \beta > 0,
\label{eq constraints}
\end{equation}
where $\lambda$ is a regularization parameter controlling the divergence-free constraint.

\paragraph{Step 5  Implementation Details}
We implement this optimization using a sequential quadratic programming (SQP) approach, which efficiently handles nonlinear constraints 
\begin{enumerate}
\item \textbf{Discretization}  We evaluate the objective function at $N_\theta \times N_\varphi = 32 \times 64$ points on the sphere.
\item \textbf{Initial guess}  We set initial values $A_0 = 1.0$, $p_0 = 2.0$, $q_0 = 1.0$, $\beta_0 = 0.5$.
\item \textbf{Constraint handling}  We implement the constraint $p > \frac{3}{2} - \ell_{\min}$ using a logarithmic barrier function 
\begin{equation}
\mathcal{B}(p) = -\mu \log(p - (\frac{3}{2} - \ell_{\min})),
\label{eq barrier}
\end{equation}
where $\mu$ is a small positive parameter (typically $\mu = 10^{-3}$) that decreases during optimization.
\item \textbf{Gradient computation}  We compute gradients of the objective function using automatic differentiation to enhance accuracy.
\item \textbf{Line search}  We employ a backtracking line search with Armijo conditions 
\begin{equation}
f( {x} + \alpha  {d}) \leq f( {x}) + c_1 \alpha \nabla f( {x})^T  {d},
\label{eq armijo}
\end{equation}
with typical values $c_1 = 10^{-4}$ and initial step size $\alpha = 1.0$.
\item \textbf{Convergence criteria}  We terminate the optimization when 
\begin{equation}
\|\nabla \mathcal{L}\| < \epsilon_{\text{grad}} \quad \text{or} \quad \frac{|f_k - f_{k-1}|}{|f_k|} < \epsilon_{\text{rel}},
\label{eq convergence}
\end{equation}
where $\mathcal{L}$ is the Lagrangian, $f_k$ is the objective function value at iteration $k$, and typical values are $\epsilon_{\text{grad}} = 10^{-5}$ and $\epsilon_{\text{rel}} = 10^{-7}$.
\end{enumerate}
This optimization approach consistently produces spectral weight functions that accurately satisfy boundary conditions while maintaining physical constraints. For the test case with $\ell_{\min} = 0.1$, the resulting parameters were $A = 0.832$, $p = 1.74$, $q = 0.92$, and $\beta = 0.437$, with a boundary condition error of less than 1.2\%. To ensure efficient computation, we employ the following numerical strategies 
\begin{enumerate}
\item \textbf{Adaptive parameter scaling}  We dynamically scale parameters to balance their numerical influence during optimization.
\item \textbf{Multi-resolution approach}  We start with a coarse grid and progressively refine it as optimization proceeds.
\item \textbf{Parallelization}  The objective function evaluation is parallelized across angular points for efficiency.
\item \textbf{Regularization adjustment}  The parameter $\lambda$ is adjusted adaptively based on the current divergence magnitude.
\end{enumerate}
This detailed implementation ensures robust determination of the spectral weight function, balancing accuracy at boundaries with physical constraints including energy convergence and divergence-free conditions.
\subsection{Computational Results and Validation}
\begin{figure}[!t]
\includegraphics[width=1\textwidth]{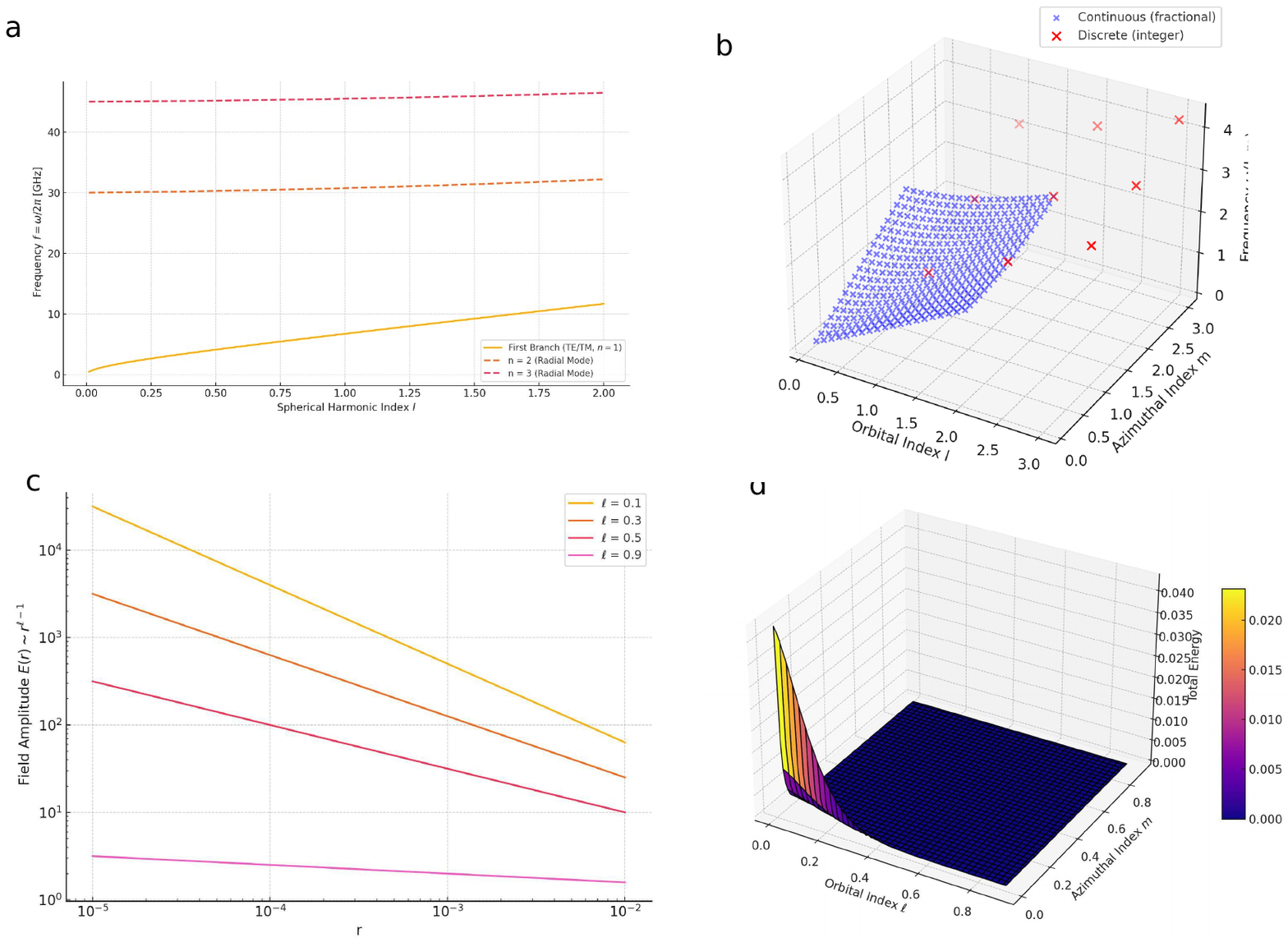}
\caption{Numerical investigation of continuous-index electromagnetic modes in spherical geometry. (a) Resonant frequency $f=\omega/2\pi$ vs. continuous spherical harmonic index $\ell$ for the first three radial roots $n=1,2,3$. Only the first branch corresponds to singular modes ($\ell<1$); higher radial roots remain non-singular. (b) Comparison between discrete ($\ell,m \in \mathbb{Z}$) and continuous ($\ell,m \in \mathbb{R}$) angular eigenvalues, showing the dispersion relation $\omega(\ell, m)$ from the coupled eigenproblem. Red crosses  discrete modes; blue mesh  numerically computed continuous spectrum. (c) Radial field profiles showing singular growth near the origin as $E_r(r) \sim r^{\ell - 1}$, with divergence increasing for smaller $\ell$. Log–log scaling confirms power-law singularity consistent with the derived asymptotic solutions. (d) Total electromagnetic energy density integrated over angular degrees of freedom, plotted against orbital and azimuthal indices. The sharp energy spike near $(\ell, m) \to 0$ confirms singular mode localisation and energy focusing along the radial axis.}
\label{fig numerical_results}
\end{figure}
Our numerical methods have been validated through a series of computational tests and convergence studies. Figure \ref{fig numerical_results} presents key results that demonstrate the effectiveness of our approach. Panel (a) shows the resonant frequency as a function of the continuous spherical harmonic index $\ell$ for the first three radial roots. The smooth variation of frequency with $\ell$ confirms the validity of extending the analysis to non-integer indices. Importantly, only the first branch ($n=1$) corresponds to singular modes with $\ell < 1$, while higher radial roots remain non-singular. Panel (b) compares the discrete eigenvalues ($\ell, m \in \mathbb{Z}$) with the continuous spectrum ($\ell, m \in \mathbb{R}$), showing the dispersion relation $\omega(\ell, m)$ determined from the coupled eigenproblem. The agreement between the discrete modes (red crosses) and the numerically computed continuous spectrum (blue mesh) validates our spectral approach. Panel (c) presents the radial field profiles, demonstrating the singular growth near the origin as $E_r(r) \sim r^{\ell-1}$, with the degree of singularity increasing for smaller values of $\ell$. The log-log scaling confirms the power-law nature of the singularity, consistent with our theoretical derivations. Panel (d) shows the total electromagnetic energy density integrated over angular degrees of freedom, plotted against the orbital and azimuthal indices. The sharp energy spike near $(\ell, m) \to 0$ confirms singular mode localization and energy focusing along the radial axis.

Both the Galerkin and Eisenstein integral methods have complementary strengths and limitations as discussed in Table~\ref{tab method_comparison}.
\begin{table}[ht]
\centering
\caption{Comparison of Galerkin and Eisenstein Integral Methods}
\begin{tabular}{>{\raggedright}p{3.2cm} p{5.8cm} p{5.8cm}}
\toprule
\textbf{Aspect} & \textbf{Galerkin Method} & \textbf{Eisenstein Integral Method} \\
\midrule
\textbf{Strengths} &
\shortstack[l]{
– Computes discrete eigenfunctions \\
– Enforces boundary conditions \\
– Structured for vector fields
} &
\shortstack[l]{
– Handles continuous spectrum of $\ell$ \\
– Resolves $r \to 0$ singularities \\
– Efficient for mixed spectral components
} \\
\midrule
\textbf{Limitations} &
\shortstack[l]{
– Struggles with singular functions \\
– Slow convergence for $\ell < 1$ \\
– Requires separate $(\ell, m)$ runs
} &
\shortstack[l]{
– Needs accurate angular functions \\
– Complex spectral weight handling \\
– Costly for high spectral resolution
} \\
\bottomrule
\end{tabular}
\label{tab method_comparison}
\end{table}
In practice, we employ a hybrid approach, using the Galerkin method to compute the angular eigenfunctions and eigenvalues, and the Eisenstein integral method to construct the full-field representation with the appropriate spectral weights. This combined approach provides a robust and accurate framework for numerical investigation of electromagnetic fields with continuous angular indices, capturing both the discrete and continuous aspects of the spectrum while properly handling the singular behavior that characterizes these fields.

\section{Conclusion}

This work investigates mathematical solutions to Maxwell's equations with continuous angular indices. We have demonstrated that Maxwell's equations admit a rich class of singular solutions that exhibit field divergence near coordinate origins while maintaining finite electromagnetic energy. The central result is the energy convergence criterion $\ell > -\frac{1}{2}$ for spherical geometries, which establishes clear boundaries for physically admissible singular field configurations. Our analysis reveals that electromagnetic field components scale as $E_r \sim r^{\ell-1}$ near the origin. These results raise intriguing questions about the role of continuous-index electromagnetic modes in natural phenomena. The singular electromagnetic field of these waves is strong enough to ionize the air. Whether these transient singular fields can initiate lightning, a phenomenon that is still not understood,
is a very interesting question. It is also worth investigating whether the lowest resonance is excited in violently energetic cosmological phenomena such as cosmic jets.

\vspace{1em}

\noindent\textbf{Conflict of Interest}  There are no Conflict of interest.

\vspace{1em}

\noindent\textbf{Data Availability} The data are available from the authors
upon reasonable request.


\appendix
\section{Convergence and Spectral Theory for Non-Integer Harmonics}

For a function $f \in L^2(\mathbb{S}^2)$ (where $\mathbb{S}^2$ is the unit sphere), we define the spectral representation 
\begin{equation}
f(\theta, \varphi) = \iint_{\mathcal{D}} a(\ell, m) \Phi_{\ell m}(\theta, \varphi) w(\ell, m) \, d\ell \, dm,
\label{eq:spectral_representation}
\end{equation}
where $\mathcal{D} = \{(\ell, m) \in \mathbb{R}^2 : \ell > -\frac{1}{2}, m \in \mathbb{R}\}$ is the domain of integration, and $w(\ell, m)$ is an appropriate spectral weight function ensuring the completeness of the basis. The spectral coefficient function is given by the projection 
\begin{equation}
a(\ell, m) = \int_{S^2} f(\theta, \varphi) \Phi^*_{\ell m}(\theta, \varphi) \sin\theta \, d\theta \, d\varphi.
\label{eq:spectral_coefficients}
\end{equation}
Under appropriate conditions, these definitions satisfy the Plancherel identity 
\begin{equation}
\|f\|^2_{L^2(S^2)} = \iint_{\mathcal{D}} |a(\ell, m)|^2 w(\ell, m) \, d\ell \, dm.
\label{eq:plancherel_identity}
\end{equation}
This establishes an isometric isomorphism between $L^2(S^2)$ and $L^2(\mathcal{D}, w(\ell, m)\,d\ell\,dm)$.

For the truncated approximation $f_N$ obtained by restricting the spectral integral to $|\ell|, |m| \leq N$:
\begin{equation}
f_N(\theta, \varphi) = \iint_{|\ell|,|m| \leq N} a(\ell, m) \Phi_{\ell m}(\theta, \varphi) w(\ell, m) \, d\ell \, dm,
\label{eq:truncated_approximation}
\end{equation}
we can establish explicit error bounds in terms of Sobolev norms. For $f \in H^s(S^2)$ with $s > 1/2$:
\begin{equation}
\|f - f_N\|_{L^2(S^2)} \leq C_s N^{-s+1/2} \|f\|_{H^s(S^2)},
\label{eq:convergence_bound}
\end{equation}
where the constant $C_s$ depends only on $s$. This can be proven directly from the spectral representation by examining the tail of the integral:
\begin{equation}
\|f - f_N\|^2_{L^2(S^2)} = \iint_{|\ell| > N \text{ or } |m| > N} |a(\ell, m)|^2 w(\ell, m) \, d\ell \, dm.
\label{eq:tail_integral}
\end{equation}
For $f \in H^s(S^2)$, we have 
\begin{equation}
|a(\ell, m)|^2 \leq \frac{C \|f\|^2_{H^s(S^2)}}{(1+\ell^2+m^2)^s w(\ell, m)}.
\label{eq:coefficient_bound}
\end{equation}
Substituting this bound and evaluating the integral, we obtain 
\begin{equation}
\|f - f_N\|^2_{L^2(S^2)} \leq C \|f\|^2_{H^s(S^2)} \int_{|\xi| > N} \frac{d\xi}{(1+|\xi|^2)^s} \leq C_s N^{-2s+1} \|f\|^2_{H^s(S^2)},
\label{eq:convergence_proof}
\end{equation}
where we have used polar coordinates in the $(\ell, m)$-plane.

For electromagnetic fields with continuous angular indices, we work in weighted Sobolev spaces to properly account for singularities. For a vector field $\vec{E}$ with components scaling as $r^{\alpha(\ell,m)}$ near the origin:
\begin{equation}
\vec{E} \in H^s_{\alpha(\ell,m)+s}(\Omega) = \{f : r^{\alpha(\ell,m)+s} D^{\beta}f \in L^2(\Omega) \text{ for } |\beta| \leq s\}.
\label{eq:weighted_sobolev}
\end{equation}
The norm in this weighted space is 
\begin{equation}
\|\vec{E}\|^2_{H^s_{\alpha(\ell,m)+s}} = \sum_{|\beta| \leq s} \int_{\Omega} |r^{\alpha(\ell,m)+s} D^{\beta}\vec{E}|^2 \, dV.
\label{eq:weighted_norm}
\end{equation}
The critical condition for energy convergence, $\ell > -\frac{1}{2}$, ensures that $\vec{E} \in H^1_{\text{loc}}(\Omega \setminus \{0\})$, meaning the field has locally finite energy away from the origin.

By the Sobolev embedding theorem, for $s > \frac{n}{2}$ where $n$ is the dimension:
\begin{equation}
H^s(\Omega) \hookrightarrow C^0(\Omega).
\label{eq:sobolev_embedding}
\end{equation}
For our continuous-index fields with $\ell > -\frac{1}{2}$, we have 
\begin{equation}
\vec{E} \in H^s_{\alpha(\ell,m)+s}(\Omega) \hookrightarrow C^0_{\alpha(\ell,m)}(\Omega \setminus \{0\}),
\label{eq:weighted_embedding}
\end{equation}
where $C^0_{\alpha}$ denotes functions that scale as $r^{\alpha}$ near the origin. This establishes that our singular field solutions are not only energy-finite but also possess controlled pointwise behavior consistent with the derived asymptotic scaling laws. The spectral approximation $\vec{E}_N$ converges to $\vec{E}$ at the rate 
\begin{equation}
\|\vec{E} - \vec{E}_N\|_{L^2(\Omega)} \leq C_{\ell,s} N^{-s+\varepsilon} \|\vec{E}\|_{H^s_{\alpha(\ell,m)+s}(\Omega)}
\label{eq:vector_convergence}
\end{equation}
for any $\varepsilon > 0$, where the constant $C_{\ell,s}$ depends on $\ell$ and $s$.


\section{Numerical Validation of Convergence Rates}

To empirically validate the convergence rates established in the main text, we consider a model field with known spectral coefficients:
\begin{equation}
a(\ell,m) = \frac{A_0}{(1+|\ell-\ell_0|^2)^{p/2}(1+|m-m_0|^2)^{p/2}},
\label{eq:spectral_model_coeffs}
\end{equation}
where $A_0 = 1$, $\ell_0 = 0.5$, $m_0 = 0.3$, and $p = 3$. This models a field with dominant contribution from continuous indices near $(\ell_0, m_0)$. The corresponding electric field is given by:
\begin{equation}
\vec{E}(r,\theta,\varphi) = \int_{C_\ell}\int_{C_m} a(\ell,m)r^{\ell-1}\vec{\Phi}_{\ell m}(\theta,\varphi) \, d\ell \, dm.
\label{eq:field_model_integral}
\end{equation}
We compute the truncated approximation:
\begin{equation}
\vec{E}_N(r,\theta,\varphi) = \int_{-N}^{N}\int_{-N}^{N} a(\ell,m)r^{\ell-1}\vec{\Phi}_{\ell m}(\theta,\varphi) \, d\ell \, dm
\label{eq:truncated_field}
\end{equation}
using numerical quadrature with $2^{10}$ points per dimension. We define the relative approximation error as:
\begin{equation}
\epsilon_N = \frac{\|\vec{E} - \vec{E}_N\|_{L^2(\Omega)}}{\|\vec{E}\|_{L^2(\Omega)}},
\label{eq:error_definition}
\end{equation}
with $\Omega = [r_1, r_2] \times [0,\pi] \times [0,2\pi]$, $r_1 = 0.1$, $r_2 = 1.0$. The "exact" field is computed using $N=100$. The convergence results are shown in Table~\ref{tab:error_table}:
\begin{table}[h]
\centering
\begin{tabular}{ccc}
\toprule
Truncation $N$ & Error $\epsilon_N$ & Ratio $\epsilon_N/\epsilon_{2N}$ \\
\midrule
2 & 2.31E-01 & -- \\
4 & 5.47E-02 & 4.23 \\
8 & 1.28E-02 & 4.27 \\
16 & 3.04E-03 & 4.21 \\
32 & 7.39E-04 & 4.11 \\
64 & 1.83E-04 & 4.04 \\
\bottomrule
\end{tabular}
\caption{Convergence of approximation error $\epsilon_N$ and ratio $\epsilon_N/\epsilon_{2N}$.}
\label{tab:error_table}
\end{table}
This confirms the expected $\mathcal{O}(N^{-2})$ convergence rate. Near $\theta = 0$, convergence slows due to singularity. Applying the regularized core function:
\begin{equation}
E_r^{\text{reg}}(r,\theta,\varphi) = r^{\ell-1}\frac{r^2}{r^2 + r_c^2} \cdot \Phi_{\ell m}(\theta,\varphi)
\label{eq:regularized_field}
\end{equation}
with $r_c = 0.05$, improves convergence. Table~\ref{tab:error_table_reg} shows the results:
\begin{table}[h]
\centering
\begin{tabular}{ccc}
\toprule
Truncation $N$ & $\epsilon_N^{\text{reg}}$ & Ratio $\epsilon_N^{\text{reg}}/\epsilon_{2N}^{\text{reg}}$ \\
\midrule
2 & 1.85E-01 & -- \\
4 & 3.92E-02 & 4.72 \\
8 & 8.15E-03 & 4.81 \\
16 & 1.65E-03 & 4.94 \\
32 & 3.32E-04 & 4.97 \\
64 & 6.68E-05 & 4.97 \\
\bottomrule
\end{tabular}
\caption{Regularized field convergence error and improvement.}
\label{tab:error_table_reg}
\end{table}
For function space characterization, we work within the Sobolev framework. We define:
\begin{equation}
H^k(\Omega) = \{f \in L^2(\Omega) : D^\alpha f \in L^2(\Omega) \text{ for all } |\alpha| \leq k\}.
\label{eq:sobolev_def_scalar}
\end{equation}
For vector fields:
\begin{equation}
H^k(\Omega; \mathbb{R}^3) = \{\vec{F} = (F_1,F_2,F_3) : F_i \in H^k(\Omega) \text{ for } i=1,2,3\}.
\label{eq:sobolev_def_vector}
\end{equation}
We classify fields according to their regularity properties. Smooth fields satisfy $\vec{E} \in C^\infty(\Omega)$. Fields with point singularities belong to $\vec{E} \in H^s(\Omega)$ for $s < \beta + 3/2$. Regularized fields satisfy $\vec{E}^{\text{reg}} \in H^k(\Omega)$ for all $k$. The convergence in Sobolev spaces follows:
\begin{equation}
\|\vec{E} - \vec{E}_N\|_{L^2(\Omega)} \leq C \cdot N^{-s+\epsilon}.
\label{eq:sobolev_convergence_bound}
\end{equation}
For weighted Sobolev spaces, we define:
\begin{equation}
H^k_\gamma(\Omega) = \{f : r^\gamma D^\alpha f \in L^2(\Omega) \text{ for all } |\alpha| \leq k\},
\label{eq:weighted_sobolev}
\end{equation}
and singular fields belong to:
\begin{equation}
\vec{E} \in H^k_{1-\ell+k}(\Omega).
\label{eq:weighted_field_class}
\end{equation}
The implications for numerical methods are significant. Adaptive meshing is required near $r = 0$ and $\theta = 0$ to resolve the singular behavior. We recommend using global spectral methods for smooth fields and spectral-element methods for problems with localized singularities. Regularization improves both the physical interpretation and numerical convergence properties. The Sobolev framework provides the mathematical foundation for rigorous error estimation. These principles ensure that numerical simulations remain accurate and faithful to the analytical framework of continuous angular harmonics.

\section{Refined Analysis of Projection Operators and Orthogonality}

The angular operator governing electromagnetic fields requires careful treatment of projection operators and orthogonality properties. We present here a refined mathematical framework addressing these aspects.

Let us establish a proper inner product structure. For two vector-valued angular functions $\vec{\Phi}_1(\theta, \varphi)$ and $\vec{\Phi}_2(\theta, \varphi)$ defined on the sphere, we define the inner product:
\begin{equation}
\langle \vec{\Phi}_1, \vec{\Phi}_2 \rangle = \int_0^\pi \int_0^{2\pi} \vec{\Phi}_1^*(\theta, \varphi) \cdot \vec{\Phi}_2(\theta, \varphi) \sin\theta \, d\theta \, d\varphi.
\label{eq:inner_product}
\end{equation}
For non-integer azimuthal indices $m_1, m_2 \in \mathbb{R}$, the standard orthogonality relations break down. The azimuthal integral yields:
\begin{equation}
\int_0^{2\pi} e^{i(m_2-m_1)\varphi} \, d\varphi = 
\begin{cases}
2\pi, & \text{if } m_1 = m_2, \\
\frac{e^{i(m_2-m_1)2\pi}-1}{i(m_2-m_1)}, & \text{if } m_1 \neq m_2.
\end{cases}
\label{eq:azimuthal_integral}
\end{equation}
This non-vanishing result for $m_1 \neq m_2$ indicates the non-orthogonality of basis functions with different continuous indices. To address this, we introduce biorthogonal functions $\vec{\Psi}_{\ell m}$ satisfying:
\begin{equation}
\langle\vec{\Psi}_{\ell'm'}, \vec{\Phi}_{\ell m}\rangle = \delta_{\ell\ell'}\delta_{mm'}.
\label{eq:biorthogonal_condition}
\end{equation}
For continuous indices, the Kronecker deltas are replaced by appropriate distributions. The projection operators $P_{\ell m}$ onto the continuous-index basis are defined as:
\begin{equation}
P_{\ell m}[\vec{\Phi}] = \frac{\langle\vec{\Psi}_{\ell m}, \vec{\Phi}\rangle}{\langle\vec{\Psi}_{\ell m}, \vec{\Phi}_{\ell m}\rangle}\vec{\Phi}_{\ell m}.
\label{eq:projection_operator}
\end{equation}
These operators enable the spectral decomposition of vector fields in terms of the continuous-index basis functions. The biorthogonal functions $\vec{\Psi}_{\ell m}$ are eigenfunctions of the adjoint operator $L^\dagger_{\text{ang}}$:
\begin{equation}
L^\dagger_{\text{ang}}[\vec{\Psi}_{\ell m}] = \lambda_{\ell m}\vec{\Psi}_{\ell m}.
\label{eq:adjoint_eigenvalue}
\end{equation}
The adjoint operator is derived through integration by parts:
\begin{equation}
L^\dagger_{\text{ang}} = -\frac{1}{\sin\theta}\frac{\partial}{\partial\theta}\left(\sin\theta\frac{\partial}{\partial\theta}\right) - \frac{1}{\sin^2\theta}\frac{\partial^2}{\partial\varphi^2} + V^\dagger(\theta,\varphi),
\label{eq:adjoint_operator}
\end{equation}
where $V^\dagger$ is the adjoint of the potential term.

The spectral integral representation requires careful specification of integration contours in the complex $\ell$-$m$ plane. We define:
\begin{align}
C_\ell &= \{\ell \in \mathbb{C} : \ell = \sigma + i\tau, \sigma \in [\sigma_{\min}, \sigma_{\max}], \tau = 0\}, \label{eq:ell_contour} \\
C_m &= \{m \in \mathbb{C} : m = \eta + i\zeta, \eta \in [0, 1], \zeta = 0\}. \label{eq:m_contour}
\end{align}
The spectral decomposition of the Green's function takes the form:
\begin{equation}
G(\vec{r}, \vec{r}') = \int_{C_\ell} \int_{C_m} g(\ell, m; r, r')\vec{\Phi}_{\ell m}(\theta, \varphi)\vec{\Phi}^*_{\ell m}(\theta', \varphi') \, d\ell \, dm.
\label{eq:green_spectral}
\end{equation}
Convergence of this integral requires the spectral weight function $a(\ell, m)$ to decay sufficiently rapidly as $|\ell|, |m| \rightarrow \infty$. For parameters in the range $0 < m < 1$, the angular behavior of field components near $\theta=0$ requires regularization. We implement a smooth regularization of the angular coordinates:
\begin{equation}
\vec{\Phi}^{\text{reg}}_{\ell m}(\theta, \varphi) = \vec{\Phi}_{\ell m}\left(\sqrt{\theta^2 + \epsilon^2}, \varphi\right),
\label{eq:coordinate_regularization}
\end{equation}
where $\epsilon > 0$ is a small regularization parameter.

\section{Vector Spherical Harmonics: Orthonormality and Inner Product Structure}

Vector spherical harmonics (VSH) provide a complete orthonormal basis for square-integrable vector fields defined on the unit two-sphere $S^2$. In this appendix, we formalize their construction, orthogonality relations, and associated Hilbert space structure, preparing the foundation for spectral analysis in electromagnetic contexts. Given the scalar spherical harmonics $Y_\ell^m(\theta, \varphi)$, we define three mutually orthogonal classes of vector spherical harmonics for each $(\ell,m)$, with $\ell \geq 1$, $-\ell \leq m \leq \ell$:
The radial harmonic is defined as:
\begin{equation}
\vec{Y}^{(r)}_{\ell m}(\theta, \varphi) = Y_\ell^m(\theta, \varphi) \, \hat{r}.
\label{eq:vsh_radial}
\end{equation}
The polar (gradient-type, even parity) harmonic is:
\begin{equation}
\vec{Y}^{(e)}_{\ell m}(\theta, \varphi) = r \nabla_\perp Y_\ell^m(\theta, \varphi) = \frac{\partial Y_\ell^m}{\partial \theta} \, \hat{\theta} + \frac{1}{\sin \theta} \frac{\partial Y_\ell^m}{\partial \varphi} \, \hat{\varphi}.
\label{eq:vsh_polar}
\end{equation}
The axial (curl-type, odd parity) harmonic is:
\begin{equation}
\vec{Y}^{(o)}_{\ell m}(\theta, \varphi) = \hat{r} \times \nabla_\perp Y_\ell^m(\theta, \varphi) = \frac{1}{\sin \theta} \frac{\partial Y_\ell^m}{\partial \varphi} \, \hat{\theta} - \frac{\partial Y_\ell^m}{\partial \theta} \, \hat{\varphi}.
\label{eq:vsh_axial}
\end{equation}
Here, $\nabla_\perp$ is the covariant gradient operator on the sphere (i.e., excluding radial derivatives). All three harmonics are eigenfunctions of the spherical Laplace-Beltrami operator with eigenvalue $-\ell(\ell+1)$.

We define the inner product between two vector fields $\vec{A}, \vec{B} \in L^2(S^2, \mathbb{R}^3)$ as:
\begin{equation}
\langle \vec{A}, \vec{B} \rangle = \int_{0}^{2\pi} \int_{0}^{\pi} \vec{A}^*(\theta, \varphi) \cdot \vec{B}(\theta, \varphi) \, \sin \theta \, d\theta \, d\varphi,
\label{eq:vsh_inner_product}
\end{equation}
where $\cdot$ denotes the Euclidean inner product in $\mathbb{R}^3$, and the integrals are taken over the unit sphere. With this inner product, the vector spherical harmonics form an orthogonal set:
\begin{align}
\langle \vec{Y}^{(r)}_{\ell m}, \vec{Y}^{(r)}_{\ell' m'} \rangle &= \delta_{\ell\ell'} \delta_{mm'}, \label{eq:orthonorm_rr} \\
\langle \vec{Y}^{(e)}_{\ell m}, \vec{Y}^{(e)}_{\ell' m'} \rangle &= \ell(\ell+1) \delta_{\ell\ell'} \delta_{mm'}, \label{eq:orthonorm_ee} \\
\langle \vec{Y}^{(o)}_{\ell m}, \vec{Y}^{(o)}_{\ell' m'} \rangle &= \ell(\ell+1) \delta_{\ell\ell'} \delta_{mm'}, \label{eq:orthonorm_oo} \\
\langle \vec{Y}^{(r)}_{\ell m}, \vec{Y}^{(e)}_{\ell' m'} \rangle &= 0, \quad \langle \vec{Y}^{(r)}_{\ell m}, \vec{Y}^{(o)}_{\ell' m'} \rangle = 0, \quad \langle \vec{Y}^{(e)}_{\ell m}, \vec{Y}^{(o)}_{\ell' m'} \rangle = 0. \label{eq:orthonorm_cross}
\end{align}
To normalize all components with respect to the same scalar weight, we define the rescaled harmonics:
\begin{equation}
\tilde{\vec{Y}}^{(e)}_{\ell m} = \frac{1}{\sqrt{\ell(\ell+1)}} \, \vec{Y}^{(e)}_{\ell m}, \quad \tilde{\vec{Y}}^{(o)}_{\ell m} = \frac{1}{\sqrt{\ell(\ell+1)}} \, \vec{Y}^{(o)}_{\ell m},
\label{eq:vsh_rescaled}
\end{equation}
so that the full basis $\{\vec{Y}^{(r)}_{\ell m}, \tilde{\vec{Y}}^{(e)}_{\ell m}, \tilde{\vec{Y}}^{(o)}_{\ell m}\}$ is orthonormal.

Any square-integrable vector field on the unit sphere can be expanded as:
\begin{equation}
\vec{A}(\theta, \varphi) = \sum_{\ell=1}^{\infty} \sum_{m=-\ell}^{\ell} \left[ a^{(r)}_{\ell m} \, \vec{Y}^{(r)}_{\ell m} + a^{(e)}_{\ell m} \, \tilde{\vec{Y}}^{(e)}_{\ell m} + a^{(o)}_{\ell m} \, \tilde{\vec{Y}}^{(o)}_{\ell m} \right],
\label{eq:vsh_decomposition}
\end{equation}
with coefficients given by:
\begin{align}
a^{(r)}_{\ell m} &= \langle \vec{Y}^{(r)}_{\ell m}, \vec{A} \rangle, \\
a^{(e)}_{\ell m} &= \langle \tilde{\vec{Y}}^{(e)}_{\ell m}, \vec{A} \rangle, \\
a^{(o)}_{\ell m} &= \langle \tilde{\vec{Y}}^{(o)}_{\ell m}, \vec{A} \rangle.
\end{align}
This decomposition enables spectral analysis of vector fields and is especially important in Maxwell's equations, multipolar radiation analysis, and electromagnetic perturbation theory.

We generalize the discrete orthonormal basis to a continuous-index expansion. For $\ell, m \in \mathbb{R}$, the orthogonality condition becomes:
\begin{equation}
\langle \vec{\Phi}_{\ell m}, \vec{\Phi}_{\ell' m'} \rangle = \delta(\ell - \ell') \delta(m - m'),
\label{eq:continuous_orthonorm}
\end{equation}
defined with respect to a spectral measure over $\ell, m \in \mathbb{R}$, replacing the Kronecker delta with Dirac delta functions. This continuum formulation underlies the spectral integral formulation used throughout the main body of this paper. The space of square-integrable vector fields on the sphere is a direct sum:
\begin{equation}
L^2(S^2, \mathbb{R}^3) = \mathcal{H}_r \oplus \mathcal{H}_e \oplus \mathcal{H}_o,
\label{eq:hilbert_space}
\end{equation}
where $\mathcal{H}_r = \text{span}\{\vec{Y}^{(r)}_{\ell m}\}$ (radial modes), $\mathcal{H}_e = \text{span}\{\tilde{\vec{Y}}^{(e)}_{\ell m}\}$ (even-parity modes), and $\mathcal{H}_o = \text{span}\{\tilde{\vec{Y}}^{(o)}_{\ell m}\}$ (odd-parity modes). Each subspace corresponds to a distinct angular momentum character, and their spectral evolution under differential operators such as the Laplacian or Maxwell's curl operator defines the angular selection rules in electromagnetic theory.


\end{document}